\documentclass[12pt]{article}
\usepackage{amsmath,amssymb}
\newtheorem{tm}{Theorem}[section]
\newtheorem{lm}[tm]{Lemma}
\newtheorem{co}[tm]{Corollary}
\newtheorem{re}[tm]{Remark}

\newtheorem{exm}[tm]{Example}
\newtheorem{pr}[tm]{Proposition}

\newcommand{\qed}{~~\hbox{\rule{4pt}{8pt}}}

\newcommand{\III}{{\vert \kern-.10em \vert \kern-.10em \vert}}

\makeatletter
 
 \@addtoreset{equation}{section}
\makeatother
 
\begin{document}
\setlength{\baselineskip}
{15.5pt}
\title{Strong uniqueness for both Dirichlet operators
and stochastic dynamics to Gibbs measures on a path space
with exponential interactions}
\author{
{\Large Sergio {Albeverio}}
\\
Institut f\"ur Angewandte Mathematik, HCM and SFB 611
\\
Universit\"at Bonn
\\
Endenicher Allee 60, D-53115 Bonn, Germany 
\\
e-mail: {\tt albeverio@uni-bonn.de} 
\vspace{4mm} \\
{\Large Hiroshi {Kawabi}}
\footnote{Corresponding author. }
\\
Department of Mathematics, Faculty of Science
\\
Okayama University \\
3-1-1, Tsushima-Naka, Kita-ku, Okayama 700-8530, Japan
\\
e-mail: {\tt kawabi@math.okayama-u.ac.jp} 
\vspace{2mm} \\ 
and 
\vspace{2mm} \\
{\Large Michael {R\"ockner}} 
\\
Fakult\"at f\"ur Mathematik,
Universit\"at Bielefeld \\
Universit\"atsstra{\ss}e 25, 
D-33501 Bielefeld, Germany 
\\
e-mail: {\tt roeckner@mathematik.uni-bielefeld.de}
}
\date{}
\maketitle 
\hspace{-8mm}
{\bf { Abstract}}:
We prove 
$L^{p}$-uniqueness of Dirichlet operators 
for
Gibbs measures 
on the path space $C(\mathbb R, \mathbb R^{d})$ associated 
with exponential type interactions in infinite volume
by extending an SPDE approach
presented in previous work by the last two named authors.
We also give an SPDE characterization of the
corresponding dynamics.
In particular, we prove  existence and uniqueness 
of a strong solution for the SPDE, though the self-interaction
potential is not assumed to be differentiable, hence the drift is
possibly discontinuous.
%
As examples, to which our results apply, we mention the stochastic 
quantization of
$P(\phi)_{1}$-,
${\rm exp}(\phi)_{1}$-, and trigonometric quantum fields in infinite volume. 
In particular, our results imply essential self-adjointness of the
generator of the stochastic dynamics for these models.
Finally, 
as an application of the strong uniqueness result for the SPDE,
we prove some functional inequalities 
for diffusion semigroups generated by 
the above Dirichlet operators. 
These inequalities are improvements of previous work 
by the second named author.
\vspace{2mm} \\
{\bf Mathematics Subject Classifications (2000):}~
35R15, 35R60, 46N50, 47D07  
\vspace{1mm} \\
{\bf Keywords:} Strong uniqueness, Dirichlet operator, Gibbs measure, 
${\rm exp}(\phi)_{1}$-quantum fields,
$L^{p}$-uniqueness, 
Essential self-adjointness, SPDE, 
Logarithmic Sobolev inequality.
\section{Introduction}
In recent years, there has been a growing interest in 
the study of 
infinite dimensional stochastic dynamics
associated with models of
Euclidean quantum field theory, hydrodynamics,  
and statistical mechanics,
see, e.g., 
Liskevich--R\"ockner \cite{Lisk-Rock},
Da Prato--Tubaro \cite{DT2} 
and Albeverio--Liang--Zegarli\'nski \cite{ALZ}, 
resp. Albeverio--Flandoli--Sinai \cite{AFS}, 
resp. Albeverio--Kondratiev--Kozitsky--R\"ockner \cite{AKKR-book}.
Equilibrium states of such dynamics are described by Gibbs 
measures.
The stochastic dynamics corresponding to these states is given by 
a diffusion semigroup, see, e.g., Albeverio \cite{A}.
On some minimal domain 
of smooth functions,
the infinitesimal generator of the semigroup coincides 
with the Dirichlet operator defined through
a classical Dirichlet form of gradient type with a Gibbs measure.
From an analytic point of view, it is very important to study
$L^{p}$-uniqueness of the Dirichlet operator, that is,
the question whether or not 
the Dirichlet operator restricted to the minimal domain has a 
unique closed extension in the $L^{p}$-space of the 
Gibbs measure under consideration,
which generates a $C_{0}$-semigroup.
As is well known, in the case of $p=2$,
this uniqueness is equivalent to 
essential self-adjointness.
We recall that essential self-adjointness is crucial in applications
to quantum mechanics to be sure that solutions of Schr\"odinger equations
are unique.
This kind of uniqueness problem on infinite dimensional state spaces
has been studied intensively by many authors.
In particular, we refer to the recent work \cite{KR} by the last two named
authors, where essential self-adjointness was proved
in the case of $P(\phi)_{1}$-quantum fields in infinite volume by
using an SPDE approach based on Da Prato--R\"ockner \cite{DR}.
Besides, in \cite{KR} also the relationship between the corresponding dynamics
and the $P(\phi)_{1}$-time evolution, which had been constructed as the strong solution
of a parabolic SPDE (\ref{GL}) in Iwata \cite{Iwa2}, is discussed.

The 
first objective of the present paper is to 
prove $L^{p}$-uniqueness of
the Dirichlet operator for all $p\geq1$, under much weaker conditions 
on the growth rate of 
the potential function
of the Gibbs measure
by a modification of the 
SPDE approach presented in \cite{KR}.
Important new examples are
${\rm exp}(\phi)_{1}$-quantum fields
in infinite volume in the context of Euclidean quantum field theory.
These models were introduced
(for the case where $\mathbb R$ 
occurring in (\ref{formal-Gibbs}) 
below is replaced by a 2-dimensional
Euclidean space-time $\mathbb R^{2}$ and where $d=1$) 
in H\o egh-Krohn \cite{Hk}, Albeverio--H\o egh-Krohn \cite{AHk} 
and further studied e.g., in Simon \cite{Si}, Fr\"ohlich \cite{Fro}, 
Albeverio--Gallavotti--H\o egh-Krohn \cite{AGHk} and Kusuoka \cite{Ku}.
More precisely, we are concerned with Gibbs measures on an infinite 
volume path space 
$C({\mathbb R}, {\mathbb R}^{d})$ given
by the following formal expression:
\begin{eqnarray}
& &
\hspace{-25mm}
Z^{-1}\exp \Big \{-\frac{1}{2}
\int_{\mathbb R} 
\big( (-\Delta_{x}+m^{2})w(x),w(x) \big)_{\mathbb R^{d}} dx
\nonumber \\
& &
\hspace{25mm}
-\int_{\mathbb R} dx \big ( \int_{\mathbb R^{d}} 
e^{(w(x), \xi )_{\mathbb R^{d}}} 
\nu(d\xi) \big ) 
\Big \} \prod_{x\in \mathbb R} dw(x).
\label{formal-Gibbs}
\end{eqnarray}
Here $Z$ is a normalizing constant, 
$m>0$ denotes mass,
$\Delta_{x}:=d^{2}/dx^{2}$, $\nu$ is a bounded positive measure on $\mathbb R^{d}$ 
with compact support, and $\prod_{x\in \mathbb R} dw(x)$ stands for a (heuristic)
volume measure on the space of maps from $\mathbb R$ into $\mathbb R^{d}$.
This has the interpretation of a quantized $d$-dimensional vector field
with an interaction of exponential type in the $1$-dimensional space-time
$\mathbb R$,
a model which is known as stochastic quantization of the ${\rm exp}(\phi)_{1}$-quantum
field model (with weight measure $\nu$).
We should mention that essential self-adjointness of the Dirichlet 
operators for such ${\rm exp}(\phi)_{1}$-quantum fields was not known yet, 
although the corresponding stochastic dynamics was constructed 
by using the Dirichlet form theory in Albeverio--R\"ockner \cite{AR}
(see also Hida--Kuo--Potthoff--Streit \cite{HKPS} for an approach based on white noise calculus).
Another important new example we handle is the model of trigonometric interactions, defined 
analogously to
(\ref{formal-Gibbs}), but with $e^{(w(x), \xi )_{\mathbb R^{d}}} $ replaced by
$\cos \{(w(x), \xi )_{\mathbb R^{d}}+\alpha \}, \alpha \in \mathbb R$.
Such a model was studied (with $\mathbb R$ replaced by a 2-dimensional space-time $\mathbb R^{2}$
and assuming $d=1$) e.g., in Albeverio--H\o egh-Krohn \cite{AHk73}, 
Fr\"ohlich 
\cite{Fro} and Albeverio--Haba--Russo \cite{AHR}.
In the present paper, we, in particular, prove essential self-adjointness of the corresponding
Dirichlet operator for all these models. 
As a consequence, the Dirichlet operator
associated with the superposition of polynomial, exponential and trigonometric
interactions, is also essentially self-adjoint. 
%

The second objective of the present paper is to discuss a
characterization of the stochastic dynamics corresponding to the above Dirichlet operator. 
%
Due to general theory, the stochastic dynamics constructed through the 
Dirichlet form approach solves the parabolic SPDE (\ref{GL}) weakly. 
However, we prove something much better, namely existence and
uniqueness of a strong solution. We achieve this by first proving
pathwise uniqueness for SPDE (\ref{GL}) and then applying the
recent work of Ondrej\'at \cite{Ondre} 
on the Yamada--Watanabe theorem for mild solutions of 
SPDE.
%
As a consequence, we have the existence of a unique strong solution to SPDE (\ref{GL})
by using simple and straightforward arguments which do not rely on any finite volume 
approximations discussed in \cite{Iwa2} in case of polynomial (i.e., smooth)
self-interaction.
%

Here we would like to emphasize that neither of the two uniqueness statements in
Theorems \ref{ES} and \ref{ES2} respectively implies the other
(cf. Remark \ref{compare} below).
%

The organization of this paper is as follows:
In Section 2, we present the framework and state our results.
There, we construct Gibbs measures as (\ref{formal-Gibbs}) rigorously
by using $d$-dimensional Brownian motion and the ground states of 
Schr\"odinger operators on $L^{2}(\mathbb R^{d}, \mathbb R)$.
After introducing our Dirichlet form and the corresponding Dirichlet operator, 
we state our main results (Theorems \ref{ES} and \ref{ES2}).
Section 3 contains the proofs, in which, we prove our main theorems. 
In our proof, we regard the Dirichlet operator as a perturbation
of the infinite dimensional Ornstein--Uhlenbeck operator by a 
possibly discontinuous and unbounded drift term.
Then we implement a modification
of a technique developed in {\cite{KR}} which in turn is based on
beautiful results of Da Prato, Tubaro and Priola in \cite{D, DT, priola}
for Lipschitz perturbations of
the Ornstein--Uhlenbeck operators.
(For other works on perturbed infinite dimensional Ornstein--Uhlenbeck
operators, see also, e.g., Albeverio--R\"ockle--Steblovskaya \cite{ARS} and 
references therein.)
To handle our quite singular drift term, the first thing to do is 
to check its $L^{p}$-integrability.
For this purpose, we make use of the asymptotic behavior for the ground state of 
the Schr\"odinger operator at infinity which, through the Feynman--Kac formula, 
has a close connection with
the growth rate of the potential function. 
We introduce an approximation scheme for the potential function by combining
the Moreau--Yosida approximation (\ref{Moreau-Yosida}) with a further 
regularization (\ref{further})
inspired by \cite{DR, KR},
and this scheme works efficiently in our proof.
To show existence and uniqueness of a strong solution to SPDE (\ref{GL}), we firstly
identify our diffusion process 
as a weak solution to an infinite system of SDEs. Secondly, we translate the
infinite dimensional SDE into the weak form of SPDE (\ref{GL}), and show 
pathwise uniqueness for it based on Marinelli--R\"ockner \cite{Carlo}.
In Section 4, we discuss some functional inequalities including the
logarithmic Sobolev inequality (\ref{LSI-ineq}) as an application of Theorem \ref{ES2},
and in Section 5, we give another proof of the
logarithmic Sobolev inequality (\ref{LSI-ineq}) by using
Lemmas \ref{Hirokawa} and \ref{Conv-GS}
on the approximation of the ground state. These lemmas play key roles when we combine 
some tightness arguments with the previous work to derive 
inequality (\ref{LSI-ineq}).
\section{Framework and Results}
We begin by introducing some notation and objects we will be working with.
We define a weight function $\rho_{r}\in C^{\infty}({\mathbb R}, 
{\mathbb R} ), r\in {\mathbb R}$, by 
$\rho_{r}(x):=e^{r\chi (x)}$, $x\in {\mathbb R}$,  
where 
$\chi \in C^{\infty}({\mathbb R}, {\mathbb R})$ 
is a positive symmetric convex function
satisfying $\chi (x)=\vert x \vert$ for $\vert x \vert \geq 1$.
We fix a positive constant $r$ sufficiently small. In particular,
we take $r>0$ such that $2r^{2}<K_{1}$ if $K_{1}>0$,
where the constant $K_{1}$ appears
in condition {\bf (U1)} below.
We set $E:=L^{2}({\mathbb R}, {\mathbb R}^{d};\rho_{-2r}(x)dx)$. This space
is a Hilbert space with its inner product defined by
$$ (w, {\tilde w})_{E}:=\int_{\mathbb R}
\big( w(x),{\tilde w}(x) {\big )}_{{\mathbb R}^{d}}
\rho_{-2r}(x)dx,\quad  w, {\tilde w}
\in E.
$$
Moreover, we set $H:=L^{2}({\mathbb R},{\mathbb R}^{d})$
and denote by $\Vert \cdot \Vert_{E}$ and $\Vert \cdot \Vert_{H}$ 
the corresponding norms in $E$ and $H$, respectively.  
We regard the dual space $E^{*}$ of $E$ as $L^{2}({\mathbb R}, {\mathbb R}^{d};\rho_{2r}(x)dx)$.
We endow $C({\mathbb R}, \mathbb R^{d})$ with the compact uniform topology and 
introduce a tempered subspace 
\begin{equation}
{\cal C}:=\{ w \in C({\mathbb R}, \mathbb R^{d} )
 \vert~\lim_{\vert x \vert \to \infty} \vert w(x) \vert \rho_{-r}(x)<\infty 
\mbox{ for every }
r>0 \}.
\label{temper}
\nonumber
\end{equation}
We easily see that the inclusion 
${\cal C} \subset 
E \cap C({\mathbb R}, \mathbb R^{d})$ is dense with respect to the topology of $E$.
Let $\cal B$ be the topological $\sigma$-field on $C({\mathbb R},\mathbb R^{d})$. 
For $T_{1}<T_{2}\in \mathbb R$,
we define by
${\cal B}_{[T_{1},T_{2}]}$ and 
${\cal B}_{[T_{1},T_{2}],c}$ the sub-$\sigma$-fields of $\cal B$
generated by 
$\{ w(x);T_{1}\leq x \leq T_{2} \}$ and $\{w(x); x\leq T_{1}, x\geq T_{2} \}$, respectively.
For $T_{1}, T_{2} \in \mathbb R$ and $z_{1},z_{2} \in {\mathbb R}^{d}$, let 
${\cal W}_{[T_{1},T_{2}]}^{z_{1},z_{2}}$ be 
the path space measure of the Brownian bridge such that $w(T_{1})=z_{1}, w(T_{2})=z_{2}$.
We sometimes regard this measure as a probability measure on the measurable space
$(C({\mathbb R}, {\mathbb R}^{d}), {\cal B})$ by considering $w(x)=z_{1}$ for
$x\leq T_{1}$ and $w(x)=z_{2}$ for $x\geq T_{2}$.
%

Following Simon \cite{simon} and Iwata \cite{Iwa1},
we now proceed to introduce rigorously
the Gibbs measure on $C({\mathbb R}, {\mathbb R}^{d})$. 
In this paper, we impose 
the following conditions on the potential function $U \in C(\mathbb R^{d}, 
\mathbb R)$:
\vspace{2mm} \\
{\bf (U1)}\quad 
There exist a constant $K_{1}\in {\mathbb R}$ and a convex function
$V:{\mathbb R}^{d} \to \mathbb R$ such that
$$ U(z)=\frac{K_{1}}{2}\vert z \vert^{2}+V(z), \qquad z\in 
\mathbb R^{d}.
$$
\vspace{-4mm} \\
{\bf (U2)} \quad There exist $K_{2}>0$, $R>0$ and $\alpha>0$ such that
$$ U(z)\geq K_{2}\vert z \vert^{\alpha}, \qquad \vert z \vert > R.
$$
\vspace{-4mm} \\
{\bf (U3)}\quad  There exist $K_{3}, K_{4}>0$ and $0<\beta <1+\frac{\alpha}{2}$ such that
$$
\vert {\widetilde {\nabla}} U(z) \vert \leq K_{3}
\exp (K_{4}\vert z \vert^{\beta}), \qquad z \in {\mathbb R}^{d},$$
where ${\widetilde {\nabla}} U(z):=K_{1}z+\partial_{0}V(z)
, z\in \mathbb R^{d}$ and $\partial_{0}V$ is the minimal section
of the subdifferential $\partial V$.
(The reader is referred to Showalter \cite{Show} for 
the definition
of the subdifferential for a convex function and its minimal section.
In the case where $U\in C^{1}(\mathbb R^{d}, \mathbb R)$,
${\widetilde {\nabla}}U$ coincides with 
the usual gradient $\nabla U$.)

Let $H_{U}:=-\frac{1}{2}\Delta_{z}+U$ be the Schr\"odinger operator on 
$L^{2}({\mathbb R}^{d}, {\mathbb R})$, where 
$\Delta_{z}:=
\sum_{i=1}^{d}{\partial^{2}}/{\partial z_{i}^{2}}$ 
is the 
$d$-dimensional Laplacian. 
Then condition {\bf (U2)} assures that
$H_{U}$ has purely discrete spectrum and a complete 
set of eigenfunctions (see, e.g., Reed--Simon \cite{rs}). 
We denote by $\lambda_{0}(>\min U)$ the minimal eigenvalue
and by $\Omega$ the corresponding normalized eigenfunction in 
$L^{2}({\mathbb R}^{d}, {\mathbb R})$. This eigenfunction is called ground state and 
it can be chosen to be strictly positive. Moreover, it has exponential decay at infinity.
To be precise, there exist some positive constants $D_{1}, D_{2}$ such that
\begin{equation}
0< \Omega(z)  \leq D_{1} \exp \big( -D_{2}\vert z\vert \hspace{0.5mm} U_{\frac{1}{2} \vert z \vert}(z)^{1/2} \big),
\quad z\in \mathbb R^{d},
\label{falloff}
\end{equation}
where $U_{\frac{1}{2} \vert z \vert}(z):=
\inf \{ U(y)\vert~
\vert y-z \vert \leq \frac{1}{2}\vert z\vert 
\}$.
See \cite[Corollary 25.13]{simon} for details. 
%

\vspace{2mm}
We are going to define a probability measure $\mu$ on $(C({\mathbb R}, {\mathbb R}^{d}), {\cal B})$.
For $T_{1}<T_{2}$, and for all
$T_{1}\leq x_{1}<x_{2}<\cdots <x_{n}\leq T_{2},~A_{1},A_{2},\cdots,A_{n} \in {\cal B}({\mathbb R}^{d})$,
we define a cylinder set $A\in {\cal B}_{[T_{1},T_{2}]}$ by
$A:=\{w\in C(\mathbb R, \mathbb R^{d})~\vert~w(x_{1})\in A_{1}, w(x_{2})\in A_{2}, \cdots, w(x_{n})\in A_{n} \}$.
Next, we set
\begin{eqnarray}
\mu(A)&:=&
\Big(
\Omega,
 e^{-(x_{1}-T_{1})(H_{U}-\lambda_{0})}
%
%
\big ({ \bf 1}_{A_{1}}
 e^{-(x_{2}-x_{1})(H_{U}-\lambda_{0})} 
\big (
{\bf 1}_{A_{2}}
 \cdots
\nonumber \\
&\mbox{ }&
\hspace{30mm}
e^{-(x_{n}-x_{n-1})(H_{U}-\lambda_{0})} 
\big ({\bf 1}_{A_{n}}
%
%
 e^{-(T_{2}-x_{n})(H_{U}-\lambda_{0})}
\Omega
\big )
\big )
\big)
\Big)_{L^{2}({\mathbb R}^{d},{\mathbb R})}
\nonumber \\
&=&
e^{\lambda_{0}(T_{2}-T_{1})}\int_{{\mathbb R}^{d}} dz_{1} \Omega(z_{1})
\int_{{\mathbb R}^{d}} dz_{2} \Omega(z_{2})p(T_{2}-T_{1},z_{1},z_{2})
\nonumber \\
&\mbox{ }& \times
\int_{C(\mathbb R, \mathbb R^{d})}
{\bf 1}_{A}(w)  \exp \big( -\int_{T_{1}}^{T_{2}}U(w(x))dx  \big) {\cal W}_{[T_{1},T_{2}]}^{z_{1},z_{2}}(dw),
\label{Construction-Gibbs}
\end{eqnarray}
where $p(t,z_{1},z_{2}), t>0, z_{1},z_{2}\in \mathbb R^{d}$,
is the transition probability density of standard Brownian motion on ${\mathbb R}^{d}$,
and we used the Feynman--Kac formula for the second line.
Then by recalling that $e^{-tH_{U}}\Omega=e^{-t\lambda_{0}}\Omega,
\Vert \Omega \Vert_{L^{2}(\mathbb R^{d},\mathbb R)}=1$ and by
the Markov property of the $d$-dimensional Brownian motion,
(\ref{Construction-Gibbs}) defines a consistent family of probability measures, 
and hence $\mu$ can be extended to a probability measure on $C({\mathbb R}, {\mathbb R}^{d})$.

In the same way as \cite[Proposition 2.7]{Iwa1}, 
we can see that
$\mu({\cal C})=1$ and 
the following DLR-equations hold
even if we replace the potential function with polynomial growth by the one 
satisfying the much weaker condition {\bf (U3)}:
\begin{eqnarray}
{\mathbb E}^{\mu} \big [ {\bf 1}_{A} \vert {\cal B}_{[T_{1},T_{2}],c} \big ](\xi)
\hspace{-1.5mm}
&=&
\hspace{-1.5mm}
Z^{-1}_{[T_{1},T_{2}]}(\xi)
\int_{A} \exp \Big( -\int_{T_{1}}^{T_{2}}U(w(x))dx  \Big)
{\cal W}_{[T_{1},T_{2}]}^{\xi(T_{1}),\xi(T_{2})}(dw),
\nonumber \\
&\mbox{ }& 
\mu \mbox{-a.e. }\xi
\in C(\mathbb R, \mathbb R^{d}), \mbox{ for all }A\in {\cal B}_{[T_{1},T_{2}]}, T_{1}<T_{2},
\label{DLR}
\end{eqnarray}
where
$Z_{[T_{1},T_{2}]}(\xi):={\mathbb E}^{{\cal W}_{[T_{1},T_{2}]}^{\xi(T_{1}),\xi(T_{2})}
}
[\exp(-\int_{T_{1}}^{T_{2}} U(w(x))dx) ]
$
is a normalizing constant.
By the continuity of 
the inclusion map of $\cal C$ into $E$, we may regard
$\mu$ as a probability measure on $E$ by identifying it with its
image measure under the inclusion map, and using that,
${\cal C} \in {\cal B}(E)$ and 
${\cal B}(E) \cap {\cal C} ={\cal B}({\cal C})$ by Kuratowski's theorem.
%
The DLR-equations (\ref{DLR}) imply that the Gibbs measure $\mu$
is 
$C^{\infty}_{0}({\mathbb R}, {\mathbb R}^{d})$-quasi-invariant, i.e.,
$\mu(\cdot +k)$ and $\mu$ are mutually equivalent, and
$
\mu(k+dw)=
\Lambda(k,w)
\mu(dw)
$
holds for every $k\in C^{\infty}_{0}({\mathbb R}, {\mathbb R}^{d})$. 
In particular by Albeverio--R\"ockner \cite[Proposition 2.7]{AR1}, $\mu(O)>0$ for every open $\emptyset \neq O \subset E$,
i.e., the topological support $\rm{supp }(\mu)$ is equal to all of $E$. 
The Radon-Nikodym density 
$\Lambda (k,w)$ is represented by
\begin{eqnarray}
\Lambda(k,w)&=&
 \exp \Big \{
\int_{\mathbb R} \Big (
U\big(w(x)\big)-U\big(w(x)+k(x)\big)
\nonumber \\
&\mbox{ }& 
\qquad \qquad \qquad 
\qquad
-\frac{1}{2}
\big \vert \frac{dk}{dx}(x) \big \vert^{2}+
(w(x), \Delta_{x}k(x)
)_{{\mathbb R}^{d}}
\Big )dx \Big \}.
 \label{quasi} 
\end{eqnarray}
%

We give the following examples which satisfy our conditions {\bf (U1)},
{\bf (U2)} and {\bf (U3)}.
\begin{exm} [$P(\phi)_{1}$-quantum fields]
\label{poly-potential case}
We consider the case where the potential function $U$ is written as 
the following potential function on $\mathbb R^{d}$:
$$ U(z)=\sum_{j=0}^{2n}a_{j}\vert z 
\vert^{j},\quad 
a_{2n}>0,~ n\in \mathbb N. $$
Especially, in the case 
$U(z)=\frac{m^{2}}{2} \vert z \vert^{2}$, $m>0$,
the corresponding Gibbs measure $\mu$ is the Gaussian measure on ${\cal C}$ with mean $0$
and covariance operator $(-\Delta_{x}+m^{2})^{-1}$. 
It is just the (space-time) free field of mass $m$ in terms of Euclidean quantum field theory.
A double-well potential 
$U(z)=a(\vert z \vert^{4}-\vert z \vert^{2}), a>0$, is also particularly important 
from the point of view of physics. 
\end{exm}
\begin{exm} [${\rm exp}(\phi)_{1}$-quantum fields]
\label{exp-potential case}
We consider an exponential type potential function 
$U: {\mathbb R}^{d} \to {\mathbb R}$ (with weight $\nu$) given by 
\begin{equation}
U(z) =\frac{m^{2}}{2} \vert z \vert^{2}+V(z):=
\frac{m^{2}}{2} \vert z \vert^{2}
+
\int_{{\mathbb R}^{d}} 
e^{(\xi, z)_{\mathbb R^{d}}}
\nu(d\xi), \quad z\in \mathbb R^{d},
\label{U-def}
\nonumber
\end{equation}
where 
$\nu$ is a bounded positive measure with
${\rm supp}(\nu) \subset 
\{ \xi \in \mathbb R^{d} \vert~\vert \xi \vert \leq L \}$
for some $L>0$. 
We note that $U$ is a smooth strictly convex function (i.e., $\nabla^{2}U \geq m^{2}$).
Hence we can take
$K_{1}=m^{2}$,
$K_{2}=\frac{m^{2}}{2}$ and $\alpha=2$. Moreover,
\begin{equation}
\vert U(z) \vert 
\leq \frac{m^{2}}{2} \vert z \vert^{2}+ \nu({\mathbb R}^{d}) e^{L\vert z \vert}
\leq \big( \frac{m^{2}}{2L^{2}}+\nu({\mathbb R}^{d}) \big) e^{2L\vert z \vert},
\quad z\in \mathbb R^{d},
\label{U-bound}
\nonumber
\end{equation}
and
\begin{equation}
\vert \nabla U(z) \vert \leq 
m^{2}\vert z \vert +\int_{\mathbb R^{d}} \vert \xi \vert e^{(\xi,z)_{\mathbb R^{d}}} \nu(d\xi)
\leq
(\frac{m^{2}}{L}+L\nu(\mathbb R^{d}))e^{L\vert z \vert}, 
\quad
z\in \mathbb R^{d}.
\label{3-1}
\nonumber
\end{equation}
Thus we can take $\beta=1$, which satisfies $\beta <1+\frac{\alpha}{2}$ in condition {\bf (U3)}.
\end{exm}
\begin{re}
\label{Betz-Hairer}
We discuss a simple example of ${\rm exp}(\phi)_{1}$-quantum fields in the case $d=1$.
This example has been discussed in the
$2$-dimensional space-time case in {\rm{\cite{AHk}}}.
Let $\delta_{a}$ be the Dirac measure at $a\in \mathbb R$ and 
we consider 
$\nu(d\xi):=\frac{1}{2}\big ( \delta_{-a}(d\xi)
+\delta_{a}(d\xi) \big ),~a>0$. Then the corresponding potential function is $U(z)=\frac{m^{2}}{2}z^{2}
+\cosh (az)$,
and {\rm{(\ref{falloff})}} implies that
the Schr\"odinger operator $H_{U}$ has a ground state $\Omega$
satisfying
\begin{equation}
0<\Omega(z) \leq D_{1} \exp \big( -\frac{D_{2}}{\sqrt 2}
\vert z\vert \hspace{0.5mm} e^{\frac{a}{4}\vert z \vert} \big),
\quad z\in \mathbb R,
\label{BH-upper}
\end{equation}
for some $D_{1}, D_{2}>0$.
By the translation invariance of the Gibbs measure $\mu$ and {\rm{(\ref{BH-upper})}}, 
there exist positive constants $M_{1}$ and $M_{2}$ such that
\begin{eqnarray}
A_{T}&:=&
\mu \big ( \{w \in C(\mathbb R, \mathbb R) \vert~\vert w(T) \vert > \frac{4}{a} \log \log T \} \big )
\nonumber \\
&=& \int_{\vert z \vert >\frac{4}{a} \log \log T} \Omega(z)^{2}dz
\nonumber \\
&\leq&  M_{1} \exp \big \{ -
M_{2}
(\log T)(\log \log T) \big \}
=M_{1} T^{-M_{2}\log \log T}
\label{BH-upper2}
\end{eqnarray}
for $T$ large enough, and {\rm{(\ref{BH-upper2})}} implies
$ \sum_{T=1}^{\infty} A_{T}<\infty$. 
Then the first Borel--Cantelli lemma yields 
$$ \mu \big ( \{w \in C(\mathbb R, \mathbb R) \vert~\limsup_{T \to \infty} \frac{
\vert w(T) \vert}{\log \log T} \leq \frac{4}{a} \} \big )=1,
$$
and thus $\mu$ is supported by a much smaller subset of $C(\mathbb R, \mathbb R)$
than ${\cal C}$.
\end{re}
\begin{exm} [Trigonometric quantum fields]
\label{sine-potential case}
We consider a trigonometric type potential function 
$U: {\mathbb R}^{d} \to {\mathbb R}$ (with weight $\nu$) given by 
\begin{equation}
U(z) =\frac{m^{2}}{2} \vert z \vert^{2}+V(z):=
\frac{m^{2}}{2} \vert z \vert^{2}
+
\int_{{\mathbb R}^{d}} 
\cos \big \{ {(\xi, z)_{\mathbb R^{d}}}+\alpha \big \}
\nu(d\xi), \quad z\in \mathbb R^{d},
\nonumber
\end{equation}
where $\alpha\in \mathbb R$, $m>0$, and $\nu$ is a bounded signed measure with finite second absolute moment, i.e., 
$$ \vert \nu \vert (\mathbb R^{d})<\infty, \quad 
K(\nu):=\int_{\mathbb R^{d}} \vert \xi \vert^{2} \vert \nu \vert (d\xi) <\infty. $$
This potential function is smooth, and it can be regarded as a bounded perturbation of a quadratic function.
Moreover, $\nabla^{2}U \geq m^{2}-K(\nu)$ and 
$$ \vert \nabla U(z) \vert \leq m^{2}\vert z \vert +K(\nu)^{1/2} \vert \nu \vert(\mathbb R^{d})^{1/2}.
$$
This type of potential functions corresponds to quantum field models with ``trigonometric interaction"
and has been discussed especially in the {\rm{2}}-dimensional space-time case (see, e.g., 
{\rm{\cite{AHk73, Fro, HKPS}}}).
\end{exm}
\begin{re} Further examples can be obtained by considering $U:\mathbb R^{d} \to \mathbb R$
of the form $U(z)=\lambda_{1}U_{1}(z)+\lambda_{2}U_{2}(z)+\lambda_{3}U_{3}(z)$, where 
$\lambda_{i}\geq 0,~i=1,2,3$, and $U_{1}$, resp. $U_{2}$, resp. $U_{3}$, is as given in
Example {\rm{\ref{poly-potential case}}}, resp. Example {\rm{\ref{exp-potential case}}}, 
resp. Example {\rm{\ref{sine-potential case}}}.
\end{re}

Now we are in a position to introduce the
pre-Dirichlet form $({\cal E},{\cal FC}_{b}^{\infty})$.
Let $K\subset E^{*}$ be a dense linear subspace of $E$ and
let ${\cal FC}_{b}^{\infty}(K)$ be the space of all smooth cylinder functions
on $E$ having the form
$$
F(w)=f(\langle w,\varphi_{1} \rangle , \ldots , \langle w,\varphi_{n} \rangle), \quad w\in E,
$$
with $n\in {\mathbb N}$, $f
\in C^{\infty}_{b}({\mathbb R}^{n}, {\mathbb R})$ 
and $\varphi_{1}, \ldots , \varphi_{n} \in K$.
Here we set $\langle w,\varphi \rangle:=\int_{\mathbb R}
(w(x), \varphi(x))_{\mathbb R^{d}}dx$ if the integral converges absolutely,
and set ${\cal FC}_{b}^{\infty}:={\cal FC}_{b}^{\infty}(C^{\infty}_{0}({\mathbb R}, {\mathbb R}^{d}))$ 
for simplicity. Since we have supp$(\mu)=E$, two different functions in
${\cal FC}_{b}^{\infty}(K)$ represent two different $\mu$-classes.
Note that ${\cal FC}_{b}^{\infty}$ is dense in $L^{p}(\mu)$ for all $p\geq 1$.
For $F \in {\cal FC}_{b}^{\infty}$, we define the $H$-Fr\'echet derivative 
$D_{H}F:E\to H$ by
$$
D_{H}F(w):=\sum_{j=1}^{n}\frac{\partial {f}}{\partial \alpha_{j}} 
(\langle w,\varphi_{1} \rangle , \ldots,
\langle w,\varphi_{n} \rangle)\varphi_{j}.
$$
Then we consider the pre-Dirichlet form $({\cal E},{\cal FC}_{b}^{\infty})$
which is given by 
$$
{\cal E}(F,G)=
\frac{1}{2} 
\int_{E} \big( D_{H}F(w), D_{H}G(w) \big)_{H}
 \mu (dw),~~F,G\in {\cal FC}_{b}^{\infty}.
$$
%
\begin{pr} \label{IbP}
\begin{equation}
{\cal E}(F,G)=-\int_{E} {\cal L}_{0}F(w) G(w) \mu(dw),
\quad F,G\in {\cal FC}_{b}^{\infty},
\label{IbP3}
\end{equation}
where ${\cal L}_{0}F\in L^{p}(\mu),~ p\geq 1,~F\in {\cal FC}^{\infty}_{b}$, is given by
\begin{eqnarray}
\hspace{-5mm}
{\cal L}_{0}F(w)
&=&\frac{1}{2} {\rm Tr}(D_{H}^{2}F(w))
+\frac{1}{2}\big \langle w, \Delta_{x}D_{H}F(w(\cdot)) \big \rangle
-\frac{1}{2}\big \langle ({\widetilde {\nabla}} U)(w(\cdot)), D_{H}F(w) 
\big \rangle
\nonumber \\
&=&
\frac{1}{2} 
\sum_{i,j=1}^{n}\frac{\partial^{2}f}
{\partial \alpha_{i}\partial \alpha_{j}}
\big(\langle w,\varphi_{1} \rangle , \ldots,
\langle w,\varphi_{n} \rangle \big )\langle \varphi_{i}, \varphi_{j}\rangle
\nonumber \\
&\mbox{ }&+
\frac{1}{2} 
\sum_{i=1}^{n}
\frac{\partial f}
{\partial \alpha_{i}}
\big (\langle w,\varphi_{1} \rangle , \ldots,
\langle w,\varphi_{n} \rangle \big )
\cdot
\big\{\langle w, \Delta_{x}\varphi_{i}\rangle -\langle 
({\widetilde {\nabla}}U)(w(\cdot)), \varphi_{i}\rangle \big\}.
\nonumber 
\end{eqnarray}
for $F(w)=f(\langle w,\varphi_{1} \rangle , \ldots,
\langle w,\varphi_{n} \rangle)$.
\end{pr}

This proposition means that the operator ${\cal L}_{0}$ is the pre-Dirichlet operator
which is associated with the pre-Dirichlet form $({\cal E},{\cal FC}_{b}^{\infty})$.
In particular, $({\cal E},{\cal FC}_{b}^{\infty})$ is closable in $L^{2}(\mu)$.
Let us denote by ${\cal D(E)}$ the completion of ${\cal FC}_{b}^{\infty}$
with respect to the ${\cal E}_{1}^{1/2}$-norm. 
(Here we use standard notations of the 
theory of Dirichlet forms, see, e.g., \cite{A, FOT, MR}.)
By standard theory (cf. \cite{A, AR90, FOT, MR}),
$({\cal E}, {\cal D(E)})$ 
is a Dirichlet form and 
the operator ${\cal L}_{0}$ has a self-adjoint extension
$({\cal L}_{\mu}, {\rm Dom}({\cal L}_{\mu}))$, called the
Friedrichs extension, corresponding to the 
Dirichlet form $({\cal E}, {\cal D(E)})$.
The semigroup $\{ e^{t{{\cal L}_{\mu}}} \}_{t\geq 0}$ generated by 
$({\cal L}_{\mu}, {\rm Dom}({\cal L}_{\mu}))$ is Markovian, i.e.,
$0\leq e^{t{{\cal L}_{\mu}}}F \leq 1$, $\mu$-a.e. whenever
$0\leq F \leq 1$, $\mu$-a.e. Moreover, since 
$\{ e^{t{{\cal L}_{\mu}}} \}_{t\geq 0}$ is symmetric on $L^{2}(\mu)$, the Markovian
property implies that
$$ \int_{E} e^{t{{\cal L}_{\mu}}} F(w) \mu(dw) \leq \int_{E} F(w) \mu(dw), 
\quad F\in L^{2}(\mu),~F \geq 0,~ \mu \mbox{-a.e.}
$$
Hence $\Vert e^{t{{\cal L}_{\mu}}}F \Vert_{L^{1}(\mu)}
\leq \Vert F \Vert_{L^{1}(\mu)}$ holds for $F\in L^{2}(\mu)$, and 
$\{e^{t{{\cal L}_{\mu}}} \}_{t\geq 0}$
can be extended as a family of $C_{0}$-semigroup
of contractions in $L^{p}(\mu)$ for all $p\geq 1$. 
See e.g., Shigekawa \cite[Proposition 2.2]{shige} for details.

On the other hand, it is a fundamental question whether the
Friedrichs extension is the only closed extension
generating a $C_{0}$-semigroup on $L^{p}(\mu), p\geq 1$, which for
$p=2$ is equivalent to the fundamental problem of essential self-adjointness of 
${\cal L}_{0}$ in quantum physics (cf. Eberle \cite{Eb}). Even if $p=2$,
in general there are many 
lower bounded self-adjoint
extensions ${\widetilde {\cal L}}$ of ${\cal L}_{0}$ in 
$L^{2}(\mu)$ which therefore generate different symmetric strongly continuous
semigroups
$\{ e^{t{\widetilde {\cal L}}} \}_{t\geq 0}$. 
If, however, we have $L^{p}(\mu)$-uniqueness of ${\cal L}_{0}$ for some
$p\geq 2$, there is hence only one semigroup which is strongly continuous 
and with generator extending ${\cal L}_{0}$. 
Consequently, in this case, only one such $L^{p}$-,
hence only one such
$L^{2}$-dynamics exists, associated with the Gibbs measure $\mu$.

The following theorems are the main results of this paper. 
For the notions of ``quasi-everywhere" and ``capacity", we refer to \cite{A, FOT, MR}.
\begin{tm}
\label{ES}
{\rm{(1)}}~The pre-Dirichlet operator $({\cal L}_{0}, {\cal FC}_{b}^{\infty})$ 
is $L^{p}(\mu)$-unique for all $p \geq 1$, i.e., there exists exactly one
$C_{0}$-semigroup in $L^{p}(\mu)$ such that its generator extends 
$({\cal L}_{0}, {\cal FC}_{b}^{\infty})$.
\\
{\rm{(2)}}~There exists a diffusion process 
${\mathbb M}:=(\Theta, {\cal F}, \{ {\cal F}_{t} \}_{t\geq 0}, \{ X_{t} \}_{t\geq 0},
\{ {\mathbb P}_{w}\}_{w\in E} )$ such that
the semigroup $\{P_{t}\}_{t\geq  0}$ 
generated by 
the unique (self-adjoint) extension of $({\cal L}_{0}, {\cal FC}^{\infty}_{b})$
satisfies the following identity
for any bounded measurable function $F:E \to \mathbb R$, and $t>0$:
\begin{equation}
P_{t}F(w)=\int_{\Theta} F(X_{t}(\omega)) {\mathbb P}_{w}(d\omega), \quad \mu\mbox{-}a.s.~w\in E.
\label{suii-hangun}
\end{equation}
Moreover, 
$\mathbb M$ is the unique
diffusion process solving
the following ``componentwise" SDE:
\begin{eqnarray}
\langle X_{t}, \varphi \rangle &=& \langle w, \varphi \rangle +\langle B_{t}, \varphi \rangle +\frac{1}{2}\int_{0}^{t}
\big \{ \langle X_{s}, \Delta_{x} \varphi \rangle-\langle ({\widetilde {\nabla}}U)(X_{s}(\cdot)), \varphi \rangle \big \} ds,
\nonumber \\
&\mbox{ }&
\hspace{55mm} t>0,~\varphi \in C^{\infty}_{0}(\mathbb R, \mathbb R^{d}), ~{\mathbb P}_{w} \mbox{-a.s.},
\label{weakform}
\end{eqnarray}
for quasi-every $w\in E$
and such that its corresponding semigroup given by {\rm{(\ref{suii-hangun})}} consists of locally uniformly bounded (in $t$)
operators on $L^{p}(\mu), p\geq 1$,
where $\{B_{t}\}_{t\geq 0}$ is an $\{{\cal F}_{t}\}_{t\geq 0}$-adapted $H$-cylindrical Brownian motion
starting at zero
defined on $(\Theta, {\cal F}, \{ {\cal F}_{t} \}_{t\geq 0}, 
{\mathbb P}_{w})$
and ${\widetilde {\nabla}}U$
was defined in condition {\bf{(U3)}}.
\end{tm}
%
\begin{tm}
\label{ES2}
For quasi-every $w\in E$, the parabolic SPDE 
\begin{eqnarray}
dX_{t}(x)=\frac{1}{2} \big \{ \Delta_{x} X_{t}(x) -({\widetilde {\nabla}}U)(X_{t}(x)) \big \} dt 
+dB_{t}(x), \quad x \in \mathbb R, ~t>0,
\label{GL}
\end{eqnarray}
has a unique strong solution
$X=\{ X_{t}^{w}(\cdot)\}_{t\geq 0}$ living in $C([0,\infty),E)$.
Namely, there exists a set $S\subset E$ with ${\rm Cap}(S)=0$
such that for any $H$-cylindrical Brownian motion $\{B_{t}\}_{t\geq 0}$ starting at zero
defined on 
a filtered probability space $(\Theta, {\cal F}, \{ {\cal F}_{t} \}_{t\geq 0}, {\mathbb P})$
satisfying the usual conditions and an initial datum 
$w \in E \setminus S$,
there exists a unique $\{ {\cal F}_{t} \}_{t\geq 0}$-adapted process 
$X=\{ X_{t}^{w}(\cdot)\}_{t\geq 0}$ 
living in $C([0,\infty),E)$ satisfying {\rm{(\ref{weakform})}}.
\end{tm}
\begin{re}
\label{compare}
Obviously, the uniqueness result in Theorem {\rm{\ref{ES2}}}
implies the (thus weaker) uniqueness stated for the diffusion process $\mathbb M$
in Theorem {\rm{\ref{ES}}}.
However, it does not imply the $L^{p}(\mu)$-uniqueness of the Dirichlet
operator. This is obvious, since a priori the latter might have extensions
which generate non-Markovian semigroups which thus have no probabilistic
interpretation as transition probabilities of a process.
Therefore, neither of the two uniqueness results in Theorems {\rm{\ref{ES}}} 
and {\rm{\ref{ES2}}}, i.e., $L^{p}(\mu)$-uniqueness of the Dirichlet operator
and strong uniqueness of the corresponding SPDE respectively,
implies the other. We refer to 
Albeverio--R\"ockner {\rm{\cite[Sections 2 and 3]{Cornell}}}
and see also {\rm{\cite[Section 8]{DR}}}
for a detailed discussion.
\end{re}
\begin{re} 
If the potential function 
$U$ is a $C^{1}$-function with polynomial growth
at infinity,
Iwata {\rm{\cite{Iwa2}}} proves that
SPDE {\rm{(\ref{GL})}} 
has a unique strong solution $X^{w}=\{ X_{t}^{w}(\cdot) \}_{t\geq 0}$
living in $C([0,\infty),{\cal C})$ for every initial datum
$w\in {\cal C}$. 
On the other hand, in the case of ${\rm exp}(\phi)_{1}$-quantum fields, 
since $({\nabla} U)(w(\cdot))\notin {\cal C}$ for $w\in {\cal C}$ in general,
we cannot expect to solve SPDE {\rm (\ref{GL})} in $C([0,\infty),{\cal C})$ 
for a given initial datum $w\in {\cal C}$.
%
Hence if we 
replace the state space ${\cal C}$ by
a much smaller tempered subspace ${\cal C}_{e}$ 
such that $({\nabla} U)(w(\cdot))\in {\cal C}_{e}$ holds for $w\in {\cal C}_{e}$,
we might construct a unique strong solution to SPDE {\rm (\ref{GL})}
living in $C([0,\infty), {\cal C}_{e})$ 
for every initial datum $w\in {\cal C}_{e}$. 
%
(A possible candidate for ${\cal C}_{e}$ could be the space of all paths
behaving like
$$
\vert w(x) \vert \sim \log (\log(\log(\log(\cdots x))))
$$ at infinity.)
We will discuss this problem in the future.
\end{re}
%
%
%
\section{Proof of the Main Results}

At the beginning, we give the proof of Proposition \ref{IbP}.
\\
{\bf Proof of Proposition \ref{IbP}:} 
Firstly, we aim to prove that 
\begin{equation}
\int_{E} \Big( \int_{\mathbb R} \vert ({\widetilde{\nabla}} U)(w(x)) 
\vert^{2} \rho_{-2r}(x) dx \Big )^{p/2} \mu(dw)
<\infty,
\quad p\geq 1.
\label{p-moment}
\end{equation}
By the translation invariance of the Gibbs measure $\mu$, for every $p\geq 2$,
it holds that
\begin{eqnarray}
\lefteqn{
\int_{E} \Big( \int_{\mathbb R} \vert ({\widetilde{\nabla}} U)(w(x)) \vert^{2} \rho_{-2r}(x) dx \Big )^{p/2} \mu(dw)
}
\nonumber \\
& \leq & \int_{E} \Big \{ \Big ( \int_{\mathbb R}  \vert ({\widetilde{\nabla}} U)(w(x)) \vert^{p} \rho_{-2r}(x) dx \Big)
\big( \int_{\mathbb R} \rho_{-2r}(x) dx \big )^{\frac{p-2}{2}} \Big \}
\mu(dw)
\nonumber \\
&\leq &
\Big( \frac{1}{r} \Big )^{\frac{p-2}{2}}
\int_{\mathbb R} \Big( \int_{E} \vert ({\widetilde{\nabla}} U)(w(0)) \vert^{p} \mu(dw) \Big) \rho_{-2r}(x) dx 
\nonumber \\
& \leq & \Big (\frac{1}{r} \Big )^{p/2} \int_{\mathbb R^{d}} \vert {\widetilde{\nabla}} U(z) \vert^{p} \Omega(z)^{2}dz
\nonumber \\
&\leq &
\Big (\frac{K_{3}^{2}}{r} \Big )^{p/2} \int_{\mathbb R^{d}} \exp (pK_{4}\vert z \vert^{\beta}) \Omega(z)^{2}dz.
\label{Omega-estimate}
\end{eqnarray}
On the other hand, condition {\bf (U2)} leads to a lower bound 
$U_{\frac{1}{2} \vert z \vert}(z)\geq \frac{K_{2}}{2^{\alpha}}\vert z \vert^{\alpha}$
for $\vert z \vert \geq 2R$. Hence we can continue to bound the integral on the 
right-hand side of (\ref{Omega-estimate}) as follows:
\begin{eqnarray}
\lefteqn{
\int_{\vert z \vert \geq 2R} \exp (pK_{4}\vert z \vert^{\beta}) \Omega(z)^{2}dz
}
\nonumber \\
&\leq &
D_{1}^{2}
\int_{\vert z \vert \geq 2R} \exp \big \{ p (K_{4}\vert z \vert^{\beta}
-D_{2} \vert z \vert U_{\frac{1}{2} \vert z \vert}(z)^{1/2}  )\big \}
dz
\nonumber \\
& \leq &
D_{1}^{2}
\int_{\vert z \vert \geq 2R} \exp \Big \{ p \big(
K_{4}\vert z \vert^{\beta}- 
\frac{D_{2} K_{2}^{1/2} }{{2}^{\alpha/2} }\vert z \vert^{1+\frac{\alpha}{2}} 
\big )
\Big \}dz
< \infty,
\label{Omega-estimate-2}
\end{eqnarray}
where we used the estimate (\ref{falloff}) for the second line and 
$\beta<1+\frac{\alpha}{2}$ for the third line.

Hence by combining (\ref{Omega-estimate}) with (\ref{Omega-estimate-2}), we see 
that
the left-hand side of (\ref{Omega-estimate}) is finite for all $p\geq 2$. 
Since $\mu$ is a probability measure on $E$, we have shown that (\ref{p-moment}) 
holds for all $p\geq 1$.
In the same way, we also have
\begin{equation}
\int_{E} \Vert w \Vert_{E}^{p} \hspace{1mm} \mu(dw) <\infty, \quad p\geq 1.
\label{p-moment2}
\end{equation}

Next, we define 
$$\beta_{\varphi}(w):=\langle w, \Delta_{x}\varphi \rangle -\langle 
({\widetilde{ \nabla }}U)(w(\cdot)), \varphi \rangle,~\quad w \in E,~
\varphi \in C^{\infty}_{0}(\mathbb R, \mathbb R^{d}).$$
Then by noting that $\varphi$ has compact support, one has $\Vert \Delta_{x} \varphi \Vert_{E^{*}}
+ \Vert \varphi \Vert_{E^{*}}<\infty$, and (\ref{p-moment}) and (\ref{p-moment2}) lead to
\begin{eqnarray}
\int_{E} \vert \beta_{\varphi}(w) \vert^{p} \mu(dw)
&\leq& 2^{p-1} \big (
\Vert \Delta_{x} \varphi \Vert_{E^{*}}^{p}
+ \Vert \varphi \Vert_{E^{*}}^{p} \big )
\nonumber \\
&\mbox{ }&
\hspace{-15mm}
\times
\int_{E} \Big \{ \Vert w \Vert_{E}^{p} +
\Big( \int_{\mathbb R} \vert ({\widetilde{\nabla}} U)(w(x)) \vert^{2} \rho_{-2r}(x) dx \Big )^{p/2}
\Big \}
\mu(dw)
<\infty.
\label{beta-integrability}
\nonumber
\end{eqnarray}
Thus we have shown that ${\cal L}_{0}F \in L^{p}(\mu)$ holds for all $p\geq 1$ and 
$F\in {\cal FC}_{b}^{\infty}$. Hence the right-hand side of (\ref{IbP3}) is well-defined and
finite, and the quasi-invariance of $\mu$ yields 
(\ref{IbP3}). 
\qed
%
\vspace{2mm} 

Before proceeding to the proofs of our main theorems, we make some preparations.
We fix a positive constant $\kappa >2r^{2}$, and set 
$$ G_{t}w(x):=
\int_{\mathbb R} \frac{1}{{\sqrt{2\pi t}}}e^{-\frac{(x-y)^{2}}{2t}} w(y)dy, \quad t>0,~x\in \mathbb R.
$$
Then by \cite[Lemma 3.2]{KR}, 
we see that $\{ e^{-\kappa t/2}G_{t} \}_{t \geq 0}$ is a strongly continuous contraction
semigroup on $E$ with 
$\Vert e^{-\kappa t/2}G_{t} \Vert_{L(E,E)} \leq \exp \{-(\frac{\kappa}{2}-r^{2})t\}$.
Let $A: {\rm Dom}(A)\subset E \to E$ be the infinitesimal generator of $\{ e^{-\kappa t/2}G_{t} \}_{t \geq 0}$.
We set $e^{tA}:=e^{-\kappa t/2}G_{t}$ throughout this paper. 
By the Hille--Yosida theorem, $(A, {\rm Dom}(A))$ is $m$-dissipative and it satisfies
\begin{equation}
(Aw,w)_{E}\leq (r^{2}-\frac{\kappa}{2}) \Vert w \Vert_{E}^{2}, \quad w\in {\rm{Dom}}(A).
\label{fukuoka}
\end{equation}
\begin{lm} \label{A}
{\rm (1)} $C^{\infty}_{0}(\mathbb R, \mathbb R^{d})$ is dense in ${\rm Dom}(A)$ with respect to
the graph norm $\Vert w \Vert_{A}:=\Vert w \Vert_{E}+ \Vert Aw \Vert_{E},~w\in {\rm Dom}(A)$, 
and we have
\begin{equation}
A\varphi=\frac{1}{2}(\Delta_{x}-\kappa)\varphi, 
\quad \varphi \in C^{\infty}_{0}(\mathbb R, \mathbb R^{d}).
\label{A-hyoji}
\end{equation}
{\rm (2)} Let $A^{*}: {\rm Dom}(A^{*}) \subset E \to E$ denote the adjoint operator of 
$(A,{\rm Dom}(A))$. Then ${\rm Dom}(A^{*})={\rm Dom}(A)$. Moreover, we have 
\begin{eqnarray}
A^{*}\varphi &=& \frac{1}{2}\Delta_{x}(\rho_{-2r}\cdot \varphi ) \rho_{2r} -\frac{\kappa}{2}\varphi
\nonumber \\
&=&
A\varphi
-2r \frac{d \chi}{dx}\cdot \frac{d\varphi}{dx}+
\big \{2r^{2} \big (
\frac{d\chi}{dx} \big)^{2}-r \Delta_{x}\chi
\big \} \varphi,\quad  \varphi \in 
{C}^{\infty}_{0}(\mathbb R, \mathbb R^{d}),
\label{A*-hyoji}
\end{eqnarray}
and 
\begin{equation}
e^{tA^{*}}w(y):=e^{-\kappa t/2} \rho_{2r}(y) \cdot G_{t}(\rho_{-2r}\cdot w)(y), \quad t>0,~ y\in \mathbb R,~w\in E.
\label{eA*}
\end{equation} 
\end{lm}
{\bf Proof:}~(1) By a straightforward computation, we can easily see 
that
$C^{\infty}_{0}(\mathbb R, \mathbb R^{d}) \subset {\rm Dom}(A)$ and
that 
(\ref{A-hyoji}) holds. We introduce
$$ {\cal C}^{\infty}_{\infty}
:=\bigcap_{k=0}^{\infty}
\bigcap_{r>0} 
\Big \{\varphi \in C^{\infty}(\mathbb R, \mathbb R^{d})~\vert~
\sup_{x\in \mathbb R}
\big \vert \frac{d^{k}\varphi}{dx^{k}}(x) \big \vert \rho_{r}(x)
<\infty \Big \}.
$$
Then $C^{\infty}_{0}(\mathbb R, \mathbb R^{d}) \subset {\cal C}^{\infty}_{\infty}$
and the differential operator $A$ can be naturally extended to 
the domain ${\cal C}^{\infty}_{\infty}$ through (\ref{A-hyoji}). By using the cut-off argument 
discussed in \cite[Lemma 4.7]{KR},
we can show that $C^{\infty}_{0}(\mathbb R, \mathbb R^{d})$ is dense in ${\cal C}^{\infty}_{\infty}$
with respect to the graph norm $\Vert \cdot \Vert_{A}$.

Now, we take a function $\varphi \in {\cal C}^{\infty}_{\infty}$. Then for every $k \in \mathbb N \cup \{0 \}$ and 
$r>0$, we can find a positive constant $C(k,r)$ such that
$ \big \vert \frac{d^{k}\varphi}{dx^{k}}(x) \big \vert \leq C(k,r) \rho_{-r}(x)$ for all $x \in \mathbb R$.
Here we recall that
\begin{equation}
\int_{\mathbb R}
\frac{1}{{\sqrt{2\pi t}}}e^{-\frac{(x-y)^{2}}{2t}} 
\rho_{-2r}(y)dy\leq e^{2r^{2}t}\rho_{-2r}(x),\quad
t>0,x\in \mathbb R.
\label{heat-kihon}
\end{equation}
(cf. e.g., Da Prato--Zabcyzk \cite[Lemma 9.44]{DZ}.)
Then for every $k \in \mathbb N \cup \{0 \}$ and $r>0$, 
\begin{eqnarray}
\big \vert \frac{d^{k}}{dx^{k}}(G_{t}\varphi)(x) \big \vert \rho_{r}(x)
& \leq &
\rho_{r}(x)
\int_{\mathbb R} \frac{1}{{\sqrt{2\pi t}}}e^{-\frac{(x-y)^{2}}{2t}} 
\big \vert \frac{d^{k} \varphi}{dx^{k}}(y) \big \vert dy
\nonumber \\
& \leq & \rho_{r}(x)
\int_{\mathbb R} \frac{1}{{\sqrt{2\pi t}}}e^{-\frac{(x-y)^{2}}{2t}} \big ( C(k,2r) \rho_{-2r}(y) \big ) dy
\nonumber \\
& \leq & C(k, 2r) \rho_{r}(x) \big( e^{2r^{2}t} \rho_{-2r}(x) \big) 
\nonumber \\
& \leq & C(k, 2r) e^{2r^{2}t} <\infty, \quad x \in \mathbb R,
\label{stable-check}
\nonumber
\end{eqnarray}
where we used (\ref{heat-kihon}) for the third line and $\rho_{-2r}(x) \rho_{r}(x) = \rho_{-r}(x) \leq 1$
for the fourth line. This means that $(e^{-\kappa t/2}G_{t})({\cal C}^{\infty}_{\infty}) 
\subset {\cal C}^{\infty}_{\infty}$ 
for all $t \geq 0$, and by \cite[Theorems 1.2 and 1.3]{Eb}, we see that ${\cal C}^{\infty}_{\infty}$ is an
operator core for $A$. Hence we have shown that ${C}^{\infty}_{0}(\mathbb R, \mathbb R^{d})$ is dense in 
${\rm Dom}(A)$ with respect to the graph norm $\Vert \cdot \Vert_{A}$.
\vspace{1mm}
\\
(2)~Since (\ref{A*-hyoji}) and (\ref{eA*}) follow
by straightforward computations, it is sufficient to show the
equivalence of the graph norms $\Vert \varphi \Vert_{A}$ and 
$\Vert \varphi \Vert_{A^{*}}$ for $\varphi\in C^{\infty}_{0}(\mathbb R, \mathbb R^{d})$.

Using integration by parts,
Young's inequality $2ab \leq \delta^{-2}a^{2}+\delta^{2}b^{2}$ and that
$\Vert \frac{d\chi}{dx} \Vert_{\infty} \leq 1$, 
we obtain
\begin{eqnarray}
\big \Vert \frac{d\varphi}{dx} \big \Vert_{E}^{2} &=&-(\varphi, \Delta_{x}\varphi)_{E}
+2r\int_{\mathbb R} \big( \varphi(x), \frac{d \varphi}{dx}(x) \big)_{\mathbb R^{d}} 
\frac{d\chi}{dx}(x) \rho_{-2r}(x)dx
\nonumber \\
& \leq & \big( \frac{1}{2\delta^{2}} \Vert \varphi \Vert_{E}^{2}+\frac{\delta^{2}}{2} \Vert \Delta_{x}\varphi
\Vert_{E}^{2} \big)+ r \big( \frac{1}{2r} \big \Vert \frac{d \varphi}{dx} \big \Vert_{E}^{2}
+ 2r \Vert \varphi \Vert_{E}^{2} \big),
\nonumber 
\end{eqnarray}
which in turn implies that
\begin{equation}
\big \Vert \frac{d \varphi}{dx} \big \Vert_{E} \leq \big( 2r+\frac{1}{\delta} \big) \Vert \varphi \Vert_{E}+\delta 
\Vert \Delta_{x} \varphi \Vert_{E}, \quad \delta>0, ~
\varphi \in C^{\infty}_{0}(\mathbb R, \mathbb R^{d}).
\label{jyunbi}
\end{equation}
Recalling (\ref{A*-hyoji}), we deduce that
\begin{eqnarray}
\Vert A \varphi \Vert_{E} &\leq &
\Vert A^{*} \varphi \Vert_{E}+2r \Vert \dot \varphi \Vert_{E}
+(2r^{2}+r\Vert \Delta_{x} \chi \Vert_{\infty})\Vert \varphi \Vert_{E}
\nonumber \\
& \leq & \Vert A^{*} \varphi \Vert_{E}+2r \big \{ \big( 2r+\frac{1}{\delta} \big) \Vert \varphi \Vert_{E}+\delta 
\Vert \Delta_{x} \varphi \Vert_{E}  \big \}
+(2r^{2}+r\Vert \Delta_{x} \chi \Vert_{\infty})\Vert \varphi \Vert_{E}
\nonumber \\
& \leq & \Vert A^{*} \varphi \Vert_{E}+4r\delta \Vert A \varphi \Vert_{E}
+\big (6r^{2}+\frac{2r}{\delta}+2r\delta \kappa +r\Vert \Delta_{x} \chi \Vert_{\infty} \big) \Vert \varphi \Vert_{E},
\label{jyunbi-2}
\end{eqnarray}
where we used $\Vert \frac{d \chi}{dx} \Vert_{\infty} \leq 1$ again
for the first line and 
(\ref{jyunbi}) for the second line.
 
Now, we choose $\delta:=\frac{1}{8r}$. Then (\ref{jyunbi-2}) implies
\begin{equation}
\Vert  A \varphi \Vert_{E} \leq 2 \Vert A^{*} \varphi \Vert_{E}+ 
\big (22 r^{2}+\kappa+r\Vert \Delta_{x} \chi \Vert_{\infty} \big) \Vert \varphi \Vert_{E},
\nonumber 
\end{equation}
and by repeating a similar argument for $A^{*}\varphi$, we also have
\begin{equation}
\Vert  A^{*} \varphi \Vert_{E} \leq \frac{3}{2} \Vert A \varphi \Vert_{E}+ 
\big (22 r^{2}+\kappa+r\Vert \Delta_{x} \chi \Vert_{\infty} \big) \Vert \varphi \Vert_{E}.
\nonumber
\end{equation}
This completes the proof.
\qed
\\
{\bf Proof of Theorem \ref{ES}:} (1)~Although we mostly follow the argument in \cite{KR}, 
which in turn is based on a modification of a technique in \cite{DR}, 
we 
give an outline of the argument for the convenience of the reader.
%
We define $E_{U}:=\{ w\in E;~\Vert ({\widetilde \nabla} U)(w(\cdot)) \Vert_{E}<\infty \}$.
Then by (\ref{p-moment}), we see that $E_{U}\in {\cal B}(E)$ and $\mu(E_{U})=1$. 
We define a measurable map
${\widetilde b}:{\rm Dom}({\widetilde b}) \subset E \to E$ with
${\rm Dom}({\widetilde b})=E_{U}$ by
\begin{equation}
{\widetilde b}(w)(\cdot):=-\frac{1}{2}(\partial_{0}V)(w(\cdot))=-\frac{1}{2}\big \{ ({\widetilde \nabla} U)(w(\cdot))
-K_{1}w(\cdot) \big \}, 
\quad w\in {\rm Dom}({\widetilde b}).
\label{b-def}
\end{equation}
We note that $\mu({\rm Dom}({\widetilde b}))=1$, and 
since $V$ is convex, ${\widetilde b}$ is dissipative, i.e.,
\begin{equation}
\big ( w_{1}-w_{2}, {\widetilde b}(w_{1})-{\widetilde b}(w_{2}) \big)_{E} \leq 0, \quad w_{1},w_{2}\in 
{\rm Dom}({\widetilde b}).
\label{tilde-b-dissipative}
\end{equation}

On the other hand, we note that ${\widetilde b}$ is not continuous on $E$ in general. 
Thus we need to introduce the following regularization
scheme.
%
For $\alpha >0$, we recall the Moreau--Yosida approximation of 
$V$ which is defined by
\begin{equation} 
V_{\alpha}(z):=\inf_{y\in {\mathbb R}^{d}} \big \{ \frac{1}{2\alpha} \vert y-z \vert^{2} +V(y) \big \},
\quad z\in \mathbb R^{d}.
\label{Moreau-Yosida}
\end{equation}
Then $V_{\alpha}(z) \nearrow V(z)$ for every $z\in \mathbb R^{d}$ as $\alpha \searrow 0$.
On the other hand, $ \partial_{0}V:{\mathbb R}^{d} \to {\mathbb R}^{d}$
is maximal
dissipative by convexity of $V$.
For $\alpha >0$, we set 
$J_{\alpha}(z):=
\big(I_{\mathbb R^{d}}+\alpha \partial_{0}V \big)^{-1}(z),~z\in \mathbb R^{d}$,
and define the Yosida approximation $(\partial_{0}V)_{\alpha}:\mathbb R^{d} \to \mathbb R^{d}$ by
$$ (\partial_{0}V)_{\alpha}(z)
:=\frac{1}{\alpha}(J_{\alpha}(z)-z)
=(\partial_{0}V)(J_{\alpha}(z)),
\quad z\in \mathbb R^{d}.
$$
Then
$(\partial_{0}V)_{\alpha}$ is monotone and the following 
Lipschitz continuity holds:
$$\big \vert (\partial_{0}V)_{\alpha}(z_{1})
-(\partial_{0}V)_{\alpha}(z_{2})
\big \vert \leq 
\frac{2}{\alpha} \vert z_{1}-z_{2} \vert, \quad 
z_{1},z_{2}\in \mathbb R^{d},$$
Furthermore, it is known that $(\partial_{0}V)_{\alpha}(z)=(\partial_{0} V_{\alpha})(z), z\in \mathbb R^{d}$
(cf. e.g., \cite[Proposition 1.8]{Show}), and
\begin{eqnarray}
& &
\big \vert (\partial_{0}V)_{\alpha}(z) 
\big \vert \leq \big \vert
\partial_{0}V(z) \big \vert, \quad z\in {\mathbb R}^{d},
\label{last-1}
\\ 
\nonumber \\
& &
\lim_{\alpha \searrow 0} (\partial_{0}V)_{\alpha}(z)
=\partial_{0}V(z),
\quad z\in \mathbb R^{d}.
\label{last-2}
\end{eqnarray}
See \cite[Theorem 1.1]{Show} for details.
We define 
${\widetilde b}_{\alpha}:E\to E$ in the same way as 
${\widetilde b}$ with $\partial_{0}V$ 
replaced by $(\partial_{0}V)_{\alpha}$. Then ${\widetilde b}_{\alpha}$ 
is Lipschitz continuous and dissipative on $E$. By (\ref{last-1}) and (\ref{last-2}),
we also have
\begin{equation}
\lim_{\alpha \searrow 0}{\widetilde b}_{\alpha}(w)={\widetilde b}(w), \quad w\in {\rm Dom}({\widetilde b}).
\label{alpha-conv}
\end{equation}
However, since ${\widetilde b}_{\alpha}$
is not differentiable in general, we need to introduce a further regularization. Let 
$B:{\rm Dom}(B)\subset E \to E$ be a self-adjoint negative definite
operator such that $B^{-1}$ is of trace class.
For any $\alpha, \beta>0$, we set 
\begin{equation}
{\widetilde b}_{\alpha, \beta}(w):=\int_{E}
e^{\beta B}{\widetilde b}_{\alpha}\big( e^{\beta B}w+y\big)
N_{\frac{1}{2}B^{-1}(e^{2\beta B}-1)}(dy),\quad  w\in E,
\label{further}
\end{equation}
where $N_{Q}$ is the standard centered Gaussian measure with covariance given by
a trace class operator $Q$.
Then by applying \cite[Theorem 9.19]{DZ},
we prove that ${\widetilde b}_{\alpha, \beta}$ is dissipative, of class $C^{\infty}$,
has bounded derivatives of all orders and
\begin{equation}
\lim_{\beta \searrow 0}{\widetilde b}_{\alpha, \beta}(w)={\widetilde b}_{\alpha}(w),
\quad
\Vert {\widetilde b}_{\alpha, \beta}(w) \Vert_{E} 
\leq C_{\alpha}(1+\Vert w \Vert_{E}), 
\quad w\in E.
\label{6-00}
\end{equation}
We also define a measurable
map $b: {\rm Dom}(b) \subset E \to E$ with ${\rm Dom}(b)=E_{U}$ by
\begin{equation}
b(w):=\frac{1}{2}(\kappa -K_{1})w+{\widetilde b}(w), \quad w \in {\rm Dom}(b),
\label{def-b}
\end{equation}
and define ${b_{\alpha, \beta}}$ with 
${\widetilde b}_{\alpha, \beta}$ replacing $\widetilde b$
in (\ref{def-b}).

Now, we consider the
stochastic evolution equation on $E$ given by
\begin{eqnarray} dX_{t}&=&AX_{t}dt +b_{\alpha, \beta}(X_{t})dt+{\sqrt Q}d{W}_{t}
\nonumber \\
&=& AX_{t}dt+\frac{1}{2}(\kappa -K_{1})X_{t}dt+
{\widetilde b}_{\alpha, \beta}(X_{t})dt+{\sqrt Q}d{W}_{t}, ~t\geq 0,
\label{abstSPDE}
\end{eqnarray}
where $Q$
is a bounded linear operator on $E$ defined by
$Qw:=\rho_{-2r}\cdot w,~w\in E$, and 
$\{ W_{t} \}_{t\geq 0}$ is an $E$-cylindrical Brownian motion defined on a
fixed filtered probability space $({\Theta}, {\cal F}, \{ {\cal F}_{t} \}_{t\geq 0}, {\mathbb P})$.
Note that $Q^{-1}$ is not bounded on $E$.
This kind of equation is regarded as an abstract formulation of SPDE (\ref{GL})
in the sense of \cite{DZ}, i.e., in the mild form.
Since each $e^{tA}{\sqrt{Q}}$ is a Hilbert--Schmidt operator on $E$ and
${b}_{\alpha, \beta}$ is Lipschitz continuous on $E$, SPDE (\ref{abstSPDE}) 
has a unique mild solution $X=\{ X^{w}_{t}(\cdot) \}_{t\geq 0}$
living in $C([0,\infty),E)$ for every initial datum $w\in E$. 
Here we recall that $X$ is a mild solution to SPDE (\ref{abstSPDE}) with 
$X_{0}=w\in E$ if one has
\begin{equation} 
X_{t}=e^{tA}w+\int_{0}^{t} e^{(t-s)A} b_{\alpha, \beta}(X_{s}) ds +\int_{0}^{t} e^{(t-s)A} {\sqrt Q}dW_{s},
\quad t>0,~ {\mathbb P}\mbox{-a.s.}
\label{abstSPDE2}
\end{equation}
By a standard coupling method for SPDEs applied to (\ref{abstSPDE}), 
we see that
\begin{equation} \big \Vert X_{t}^{w}-X_{t}^{\tilde w} \big \Vert_{E}\leq
e^{\frac{(-K_{1}+2r^{2})t}{2}}
\Vert w-{\tilde w} \Vert_{E}, \qquad w, {\tilde w}\in E,
\label{coupling-est}
\end{equation} 
also holds with probability one.
We can then define the transition semigroup 
corresponding to SPDE (\ref{abstSPDE}), 
denoted by $\{P^{\alpha, \beta}_{t} \}_{t\geq 0}$.

For 
$F\in 
{\cal FC}^{\infty}_{b}$ 
and $\lambda >(-\frac{K_{1}}{2}+r^{2})\vee 0$, we consider the function
$$
\Phi_{\alpha, \beta}(w):=\int_{0}^{\infty} e^{-\lambda t} P_{t}^{\alpha, \beta}F(w) dt, \quad w\in E.
$$
Then (\ref{coupling-est}) leads us to the estimate
\begin{equation}
\Vert D\Phi_{\alpha, \beta}(w) \Vert_{E} \leq \frac{2}{2\lambda 
+K_{1}-2r^{2}} \Vert DF \Vert_{\infty}, \quad w\in E,
\label{grad-est}
\end{equation}
where $DF:E \to E$ is the $E$-Fr\'echet derivative of $F$. We have the relation 
$D_{H}F={\sqrt Q}DF$. 
By Proposition \ref{IbP}, $({\cal L}_{0}, {\cal FC}^{\infty}_{b})$ is dissipative
in $L^{p}(\mu), p\geq 1$, and then it is closable. Let $({\overline{\cal L}}_{0},
{\rm Dom}({\overline{\cal L}}_{0}) )$ denote the closure in $L^{p}(\mu)$.
However, since it is not easy to consider ${\overline{\cal L}}_{0}$ directly, we 
need to insert a tractable space between ${\cal FC}^{\infty}_{b}$ and
${\rm Dom}({\overline{\cal L}}_{0})$. 
Here we recall some beautiful results on Lipschitz perturbations of 
Ornstein--Uhlenbeck operators 
discussed in
\cite{D, DT, priola}. 
By modifying the results in \cite{D, DT, priola} for our use, we deduce that
$\Phi_{\alpha, \beta}$ belongs to a ``nice" domain ${\cal D}(L, C^{1}_{b,2}(E))$ 
(see \cite{KR} for the precise definition and details) 
of the Ornstein--Uhlenbeck operator $L$ associated with the SPDE
\begin{equation} dY_{t}=AY_{t}dt+{\sqrt Q}d{W}_{t}, ~t\geq 0.
\nonumber
\end{equation}
Moreover, recalling (\ref{p-moment2}), we see that ${\overline {\cal L}}_{0}F=LF+(b,DF)_{E}$ for
$F\in {\cal D}(L, C^{1}_{b,2}(E))$ and this identity implies the inclusion
${\cal D}(L, C^{1}_{b,2}(E)) \subset {\rm Dom}({\overline{\cal L}}_{0})$. 
Hence we have $\Phi_{\alpha, \beta} \in {\rm Dom}({\overline{\cal L}}_{0})\cap C^{2}_{b}(E)$ and 
moreover $\Phi_{\alpha, \beta}$
satisfies 
\begin{equation}
(\lambda -{\overline {\cal L}}_{0})\Phi_{\alpha, \beta}=F+
\big ({\widetilde b}_{\alpha, \beta}-{\widetilde b}, D\Phi_{\alpha, \beta} \big)_{E}.
\label{6-1}
\end{equation}
By using (\ref{grad-est}), the right-hand side of (\ref{6-1}) can be estimated as follows:
\begin{eqnarray}
I_{\alpha, \beta}&:=&\int_{E}
\big \vert
 \big ({\widetilde b}_{\alpha, \beta}(w)-{\widetilde b}(w), D\Phi_{\alpha, 
\beta}(w)  \big )_{E} \big \vert^{p} \mu(dw)  
\nonumber \\
&\leq & \Big(\frac{2}{2\lambda+K_{1}-r^{2}}\Vert DF \Vert_{\infty} \Big)^{p}
\int_{E}
\big \Vert {\widetilde b}_{\alpha, \beta}(w)-{\widetilde b}(w)
\big \Vert_{E}^{p}
\mu(dw).
\label{6-2}
\end{eqnarray}
Recalling (\ref{last-1}), (\ref{last-2}), (\ref{6-00}) and using Lebesgue's
dominated convergence theorem, we conclude that
\begin{equation}
\lim_{\alpha \searrow 0}\lim_{\beta \searrow 0}I_{\alpha, \beta}
=\lim_{\alpha \searrow 0}\big(\limsup_{\beta \searrow 0}I_{\alpha, \beta}\big)=0.
\nonumber
\end{equation}
From this and (\ref{6-1}), (\ref{6-2}), we obtain
$$ \lim_{\alpha \searrow 0}\lim_{\beta \searrow 0}
(\lambda -{\overline {\cal L}}_{0})\Phi_{\alpha, \beta}=F \quad \mbox{in }L^{p}(\mu).
$$
This means that the closure of ${\rm Range}(\lambda -{\overline {\cal L}}_{0})$
contains 
${\cal FC}^{\infty}_{b}$.
Since 
${\cal FC}^{\infty}_{b}$
is dense in $L^{p}(\mu)$, 
${\rm Range}(\lambda -{\overline {\cal L}}_{0})$ is also dense in $L^{p}(\mu)$.
Then by the Lumer--Phillips
theorem, we have that $({\overline {\cal L}}_{0}, {\rm Dom}({\overline {\cal L}}_{0}))$
generates a $C_{0}$-semigroup in $L^{p}(\mu)$, and 
this completes the proof of (1).
\vspace{2mm} \\
%
(2)~Since $C^{\infty}_{0}(\mathbb R, \mathbb R^{d})$ is dense in $E^{*}$, 
${\cal D(E)}$ coincides with the closure of ${\cal FC}^{\infty}_{b}(E^{*})$
with respect to the ${\cal E}_{1}^{1/2}$-norm. Thus, we can directly apply the 
general methods of the theory of Dirichlet forms \cite{A, MR} to prove quasi-regularity of
$({\cal E}, {\cal D(E)})$ and the existence of a diffusion process $\mathbb M$ properly
associated with $({\cal E}, {\cal D(E)})$. 

Here, following R\"ockner \cite{Ro-JFA} and Funaki \cite{Funaki}, we introduce scaled  Sobolev spaces:
$$H^{m}_{r}(\mathbb R, \mathbb R^{d}):=\{ \varphi \vert~\rho_{r}\varphi \in H^{m}(\mathbb R, \mathbb R^{d}) \},
\quad m\geq 0,~r\in \mathbb R,
$$
equipped with norms
$\vert \varphi \vert_{m,r}:=\Vert
\rho_{r}\varphi \Vert_{H^{m}(\mathbb R, \mathbb R^{d})}$.
Note that this norm is equivalent to  $\Vert \varphi \Vert_{m,r}
:=\sum_{k=0}^{m} \Vert \rho_{r} \big (\frac{d^{k} \varphi}{dx^{k}}\big) 
\Vert_{L^{2}(\mathbb R, \mathbb R^{d})}$ in the case $m\in \mathbb N \cup \{0 \}$.
Let $(H^{m}_{r}(\mathbb R, \mathbb R^{d}))^{*}$ be the dual space of $H^{m}_{r}(\mathbb R, \mathbb R^{d})$.
Then we have $$(H^{m}_{r}(\mathbb R, \mathbb R^{d}))^{*}=H^{-m}_{-r}(\mathbb R, \mathbb R^{d})
=\{ w \vert~\rho_{-r}w\in H^{-m}(\mathbb R, \mathbb R^{d}) \},$$
and, clearly $H=H^{0}_{0}(\mathbb R, \mathbb R^{d})$, $E=H^{0}_{-r}(\mathbb R, \mathbb R^{d})$.
For our later use, we consider a separable Hilbert space 
${\cal H}:=H^{-2}_{-r}(\mathbb R, \mathbb R^{d})$. Since ${\cal H}^{*}
=H^{2}_{r}(\mathbb R, \mathbb R^{d})$, we have
$$ C^{\infty}_{0}(\mathbb R, \mathbb R^{d}) \subset {\cal H}^{*}
\subset E^{*}
\subset H^{*} \equiv H \subset E \subset {\cal H}
$$
and the inclusions are dense and continuous.

Let $D:=\{\varphi_{n} \}_{n=1}^{\infty}
\subset C^{\infty}_{0}(\mathbb R, \mathbb R^{d})$ be the countable weakly dense $\mathbb Q$-linear subspace 
of ${\cal H}^{*}$ constructed on page 369 of Albeverio--R\"ockner \cite{AR}.
Then by
\cite[Theorem 5.3]{AR}, for each $n\in \mathbb N$, there exists some $S_{n} \subset E$
with ${\rm Cap}(S_{n})=0$ such that the diffusion process $\mathbb M$ satisfies
\begin{equation}
\langle X_{t}, \varphi_{n} \rangle = \langle w, \varphi_{n} \rangle +B^{(n)}_{t}
+\frac{1}{2}\int_{0}^{t}
\beta_{\varphi_{n}}(X_{s}) ds,
\quad t>0,~{\mathbb P}_{w} \mbox{-a.s.},
\label{weakform-2}
\end{equation}
for all $w \in E \setminus S_{n}$, where $\{B^{(n)}_{t} \}_{t\geq 0}$
is a one-dimensional  
$\{ {\cal F}_{t} \}$-adapted Brownian motion on $({\Theta}, {\cal F}, {\mathbb P}_{w})$
starting at zero
multiplied by $\Vert \varphi_{n} \Vert_{H}$.
On the other hand, by recalling (\ref{p-moment}), (\ref{p-moment2}) and
\cite[Lemma 4.2]{AR},
there exists a set 
$S_{0}\subset E$ with ${\rm Cap}(S_{0})=0$ such that 
\begin{equation}
{\mathbb P}_{w} \Big ( \int_{0}^{T} 
\big(
\Vert (\tilde \nabla U)(X_{s}(\cdot)) \Vert_{E} +\Vert X_{s}(\cdot) \Vert_{E}
\big)
ds <\infty \mbox{ for all } T>0 \Big )=1
\label{weakform-extend}
\end{equation}
for any $w\in E\setminus S_{0}$.
Here we set $S:=\cup_{n=0}^{\infty} S_{n}$. Obviously, ${\rm Cap}(S)=0$.
By noting that the embedding map $H \hookrightarrow {\cal H}$ is a Hilbert--Schmidt operator
(cf. \cite[Remark 2.1]{Funaki}), and \cite[Remark 6.3]{AR},
we can apply \cite[Lemma 6.1 and Theorem 6.2]{AR}, which
implies that 
there exists an $\{ {\cal F}_{t} \}_{t\geq 0}$-Brownian motion on 
$({\Theta}, {\cal F}, {\mathbb P}_{w})$ with values in
$\cal H$ 
starting at zero with covariance $(\cdot, \cdot)_{H}$ 
(i.e., an $H$-cylindrical Brownian motion)
under $\mathbb P_{w}$ for 
every $w \in E \setminus S$ such that 
\begin{equation}
\langle B_{t}, \varphi_{n} \rangle
:=
\mbox{ }_{{\cal H}} \langle B_{t}, \varphi_{n} \rangle_{{\cal H}^{*}}=B^{(n)}_{t},
\quad n \in \mathbb N,~t \geq 0,~~{\mathbb P}_{w} \mbox{-a.s.},~~w\in E\setminus S.
\label{Wiener process}
\end{equation}

Since
$D$ is dense in $C^{\infty}_{0}(\mathbb R, \mathbb R^{d})$ 
with respect to the weak topology of ${\cal H}^{*}$, for every 
$\varphi \in C^{\infty}_{0}(\mathbb R, \mathbb R^{d})$,
we can take a subsequence $\{ \varphi_{n(j)} \}_{j=1}^{\infty} \subset D$  
such that $\varphi_{n(j)} \to \varphi$ weakly in ${\cal H}^{*}$ as $j \to \infty$.
Furthermore, the Banach--Saks theorem implies that,
selecting another subsequence again denoted by $\{ \varphi_{n(j)} \}_{j=1}^{\infty}$,
the Ces\`aro mean
${\hat \varphi}_{k}:=\frac{1}{k} \sum_{j=1}^{k} \varphi_{n(j)}, k\in \mathbb N$,
converges to $\varphi$ strongly in ${\cal H}^{*}$ as $k \to \infty$. Thus
$\Vert \varphi-{\hat \varphi}_{k} \Vert_{E^{*}}+\Vert    
\Delta_{x}\varphi-\Delta_{x} {\hat \varphi}_{k} \Vert_{E^{*}} \to 0$
as $k \to \infty$. 
On the other hand, (\ref{weakform-2}) and  (\ref{Wiener process}) imply
\begin{equation}
\langle X_{t}, \hat \varphi_{k} \rangle = \langle w, \hat \varphi_{k} \rangle 
+\langle B_{t}, \hat \varphi_{k} \rangle
+\frac{1}{2}\int_{0}^{t}
\beta_{\hat \varphi_{k}}(X_{s}) ds,
\quad t>0,~{\mathbb P}_{w} \mbox{-a.s.},
\label{weakform-3}
\end{equation}
for all $w \in E \setminus S$.
Hence due to (\ref{weakform-extend}) 
we can take the limit $k \to \infty$ on both sides of
(\ref{weakform-3}) to obtain 
SDE (\ref{weakform}) for all $w\in E \setminus S$.
Besides, the uniqueness statement for $\mathbb M$
is derived from item (1) (cf. \cite[Sections 2 and 3]{Cornell} and also
\cite[Section 8]{DR}).
This completes the proof.
\qed
%
%
\vspace{2mm} \\
{\bf Proof of Theorem \ref{ES2}:}~By noting (\ref{weakform-extend}),
the fact that $Q^{-1}(C^{\infty}_{0}(\mathbb R, \mathbb R^{d}))
=C^{\infty}_{0}(\mathbb R, \mathbb R^{d})$ and Theorem \ref{ES},
we can read  
(\ref{weakform}) as 
\begin{eqnarray}
(X_{t}, \varphi)_{E} &=& (w, \varphi)_{E} +\int_{0}^{t} ({\sqrt Q}\varphi, dW_{s})_{E} +\int_{0}^{t}
\big \{ (X_{s},A^{*}\varphi)_{E}+(b(X_{s}), \varphi)_{E} \big \} ds,
\nonumber \\
&\mbox{ }&
\hspace{35mm} t>0,~\varphi \in C^{\infty}_{0}(\mathbb R, \mathbb R^{d}),
~{\mathbb P}_{w} \mbox{-a.s.},~w\in E\setminus S. 
\label{weakform-E}
\end{eqnarray}
for all $w\in E\setminus S$,
where $\{W_{t} \}_{t\geq 0}$ is an $\{ {\cal F}_{t} \}_{t\geq 0}$-adapted 
$E$-cylindrical Brownian motion corresponding to the $H$-cylindrical Brownian motion
$\{ B_{t} \}_{t\geq 0}$ defined on 
$(\Theta, {\cal F}, {\mathbb P}_{w})$. (See \cite[Remark 3.5]{KR} for details.)
Furthermore, 
by recalling Lemma \ref{A}, we have equation (\ref{weakform-E}) for every
$\varphi \in {\rm Dom}(A^{*})$.
We also mention that (\ref{weakform-E}) is equivalent to the mild-form (\ref{abstSPDE2})
of SPDE (\ref{abstSPDE}) with $b_{\alpha, \beta}$ replaced by $b$. We refer to
Ondrej\'at \cite[Theorem 13]{Ondre} for details.

Now, we prove pathwise uniqueness based on the argument of Marinelli--R\"ockner \cite{Carlo}.
Suppose that $X=X^{w}$ and ${\widetilde X}={\widetilde X}^{w}$ are two weak solutions to SPDE
(\ref{abstSPDE}) defined on the same filtered probability space 
$(\Theta, {\cal F}, \{ {\cal F}_{t} \}_{t\geq 0}, {\mathbb P})$ with the same $E$-cylindrical
Brownian motion $\{ W_{t} \}_{t \geq 0}$ and 
$X_{0}={\widetilde X}_{0}=w\in E\setminus S$
such that
\begin{equation}
\int_{0}^{T} \Vert b(X_{s}) \Vert_{E}ds <\infty, \quad  
\int_{0}^{T} \Vert b({\widetilde X}_{s}) \Vert_{E}ds <\infty \quad \mbox{ for all } T>0,~{\mathbb P}\mbox{-a.s.}
\label{p-integrability}
\end{equation}
We fix $T>0$ from now on, and set $\Psi_{t}:=X_{t}-{\widetilde X}_{t}$. Note that it enjoys an $\omega$-wise
equation
$$
d\Psi_{t}=A\Psi_{t}dt +(b(X_{t})-b({\widetilde X}_{t}))dt, \quad 0<t \leq T,
$$
with the initial datum $\Psi_{0}=0$, again to be understood in the mild form.
Since $X$ and $\widetilde X$ have continuous paths on $E$, (\ref{p-integrability}) implies that
$b(X_{\cdot})-b({\widetilde X}_{\cdot}) \in
L^{1}([0,T], E)$ and $\sup_{0\leq t \leq T} \Vert \Psi_{t} \Vert_{E} <\infty$ hold for $\mathbb P$-a.s 
$\omega \in \Theta$. 
Let $\{ \varphi_{n} \}_{n=1}^{\infty} \subset 
C^{\infty}_{0}(\mathbb R, \mathbb R^{d})$ be a CONS of $H$, 
and we set ${\widetilde \varphi}_{n}:=\rho_{r} \varphi_{n}$
and $e_{n}:=(I+\varepsilon A^{*})^{-1}{\widetilde \varphi}_{n} \in {\rm Dom}(A^{*})$ for 
$n \in \mathbb N$. We mention that $\{ {\widetilde \varphi}_{n} \}_{n=1}^{\infty}$ is a CONS of $E$.
%
%
Recalling 
(\ref{weakform-E}) and applying It\^o's formula, we have
\begin{eqnarray}
(e_{n},\Psi_{t})_{E}^{2} 
&=& 2 \int_{0}^{t} \Psi_{n}(s) d\Psi_{n}(s)
\nonumber \\
& \mbox{ }& 
+2
\int_{0}^{t} ( e_{n},  \Psi_{s} )_{E}
\big (e_{n}, b(X_{s})-b({\widetilde X}_{s}) \big )_{E} ds 
\nonumber \\
&=:&
2 \big( J_{n}^{1}(t)+J_{n}^{2}(t) \big), \qquad 0\leq t \leq T.
\label{Carlo-1}
\end{eqnarray}

For the first term $J_{n}^{1}(t)$, Lebesgue's dominated convergence theorem
leads us to
\begin{eqnarray}
\sum_{n=1}^{\infty} J_{n}^{1}(t)
&=&
\int_{0}^{t} \sum_{n=1}^{\infty}
\big( (I+\varepsilon A^{*})^{-1}{\widetilde \varphi}_{n},
\Psi_{s} \big )_{E}
\cdot
\big ( (A(I+\varepsilon A)^{-1})^{*}{\widetilde \varphi}_{n}, {\Psi}_{s} \big )_{E} ds
\nonumber \\
&=&
\int_{0}^{t} \sum_{n=1}^{\infty}
\big( {\widetilde \varphi}_{n}, 
(I+\varepsilon A)^{-1}{\Psi}_{s} \big )_{E} \cdot
\big ( {\widetilde \varphi}_{n}, A(I+\varepsilon A)^{-1}{\Psi}_{s} \big )_{E} ds
\nonumber \\
&=& \int_{0}^{t} \big( (I+\varepsilon A)^{-1}\Psi_{s},
A(I+\varepsilon A)^{-1}\Psi_{s} \big )_{E} ds
\nonumber \\
&\leq & \big (r^{2}-\frac{\kappa}{2} \big) \int_{0}^{t} 
\big \Vert (I+\varepsilon A)^{-1}\Psi_{s} \big \Vert_{E}^{2}ds,
\label{Carlo-2}
\end{eqnarray}
where we used $(I+\varepsilon A^{*})^{-1}
=( (I+\varepsilon A)^{-1})^{*}$ and the fact that $A^{*}$ and $(I+\varepsilon A^{*})^{-1}$ commute
for the second line, and (\ref{fukuoka}) for the fourth line.

For the second term $J_{n}^{2}(t)$, since we have
\begin{eqnarray}
\lefteqn{
\int_{0}^{t} \big \Vert b(X_{s})-b({\widetilde X}_{s}) \big \Vert_{E} \Vert \Psi_{s} \Vert_{E} ds 
}
\nonumber \\
&\leq &
\big( \sup_{0\leq t \leq T} \Vert  \Psi_{s} \Vert_{E} \big)
\int_{0}^{T} \big \Vert b(X_{s})-b({\widetilde X}_{s}) \big \Vert_{E} ds<\infty, \quad 0\leq t \leq T,
\quad {\mathbb P}\mbox{-a.s.},
\nonumber
\end{eqnarray}
Lebesgue's dominated convergence theorem also yields
\begin{eqnarray}
\sum_{n=1}^{\infty} J_{n}^{2}(t) 
&=&
\int_{0}^{t} \sum_{n=1}^{\infty}
\big( {\widetilde \varphi}_{n}, 
(I+\varepsilon A)^{-1}{\Psi}_{s} \big )_{E} \cdot
\big ( {\widetilde \varphi}_{n}, (I+\varepsilon A)^{-1} \big (b(X_{s})-b({\widetilde X}_{s}) \big)
\big )_{E} ds
\nonumber \\
&=&
\int_{0}^{t}
\Big(  (I+\varepsilon A)^{-1}\Psi_{s}, (I+\varepsilon A)^{-1} \big (b(X_{s})-b({\widetilde X}_{s}) \big) \Big )_{E} ds. 
\label{Carlo-3}
\end{eqnarray}
Then by putting (\ref{Carlo-2}) and (\ref{Carlo-3}) into (\ref{Carlo-1}), 
we have
\begin{eqnarray}
\Vert (I+\varepsilon A)^{-1}\Psi_{t} \Vert_{E}^{2} &=&2 \sum_{n=1}^{\infty} \big( J_{n}^{1}(t)+J_{n}^{2}(t) \big)
\nonumber \\
& \leq & (2r^{2}-{\kappa}) \int_{0}^{t} 
\big \Vert (I+\varepsilon A)^{-1}\Psi_{s} \big \Vert_{E}^{2}ds
\nonumber \\
&\mbox{  }&
+2 \int_{0}^{t}
\Big(  (I+\varepsilon A)^{-1}\Psi_{s}, (I+\varepsilon A)^{-1} \big (b(X_{s})-b({\widetilde X}_{s}) \big) \Big )_{E} ds.
\nonumber
\end{eqnarray}
Moreover letting $\varepsilon \searrow 0$ on both sides,
and recalling the dissipativity (\ref{tilde-b-dissipative}) for $\widetilde b$ and 
(\ref{def-b}), 
we obtain that
$$
\Vert \Psi_{t} \Vert_{E}^{2} \leq (-K_{1}+2r^{2}) \int_{0}^{t} 
\Vert \Psi_{s} \Vert_{E}^{2}ds.
$$
Hence, we have $\Psi_{t}=X_{t}-{\widetilde X}_{t}=0,~0\leq t \leq T$, $\mathbb P$-almost
surely
by an application of
Gronwall's inequality, which proves the pathwise uniqueness.
Then by \cite[Theorem 2]{Ondre}, a Yamada--Watanabe type argument implies that
SPDE (\ref{GL}) has a unique strong solution. This completes the proof. 
\qed
\vspace{2mm}

By repeating the same argument as in the above proof, we can easily deduce
the following coupling estimates
(\ref{flow}) and (\ref{flow-2})
which play crucial roles in the next section.
\begin{co}
\label{cor}
Let $X^{w}$ and $X^{\tilde w}$ denote the strong solutions of 
SPDE {\rm{(\ref{GL})}} with the initial
datum $X^{w}_{0}=w\in E\setminus S$ and $X^{\tilde w}_{0}=\tilde w\in E\setminus S$,
respectively. Then 
\begin{equation} \big \Vert X_{t}^{w}-X_{t}^{\tilde w} \big \Vert_{E}\leq
e^{\frac{(-K_{1}+2r^{2})t}{2}} \Vert w-\tilde w \Vert_{E}, \quad t\geq 0,~{\mathbb P}\mbox{-a.s.} 
\label{flow}
\end{equation} 
In addition, for every $h \in H\setminus S$, 
we have 
\begin{equation} \Vert X_{t}^{w+h}-X_{t}^{w} \Vert_{H}\leq
e^{-\frac{K_{1}t}{2}} \Vert h \Vert_{H}, \quad t\geq 0,~{\mathbb P}\mbox{-a.s.}
\label{flow-2}
\end{equation} 
\end{co}
\section{Some Functional Inequalities}
In this section, as an application of Theorem \ref{ES2} and 
Corollary \ref{cor},
we present some functional inequalities 
for the diffusion semigroup $\{P_{t} \}_{t\geq 0}$ generated by the Dirichlet
operator ${\cal L}_{\mu}$.  In particular, we 
prove the
gradient estimate for $\{P_{t} \}_{t\geq 0}$ and
logarithmic Sobolev inequalities under much weaker conditions
on the regularity and the growth rate of the potential function $U$ than 
in the previous papers 
\cite{kawa-POTA, Kawa-LSI}
(which however already included the $P(\phi)_{1}$-case).

There are at present some approaches to derive these functional inequalities, 
and it is well-known that
Bakry--\'Emery's $\Gamma_{2}$-method (cf. Bakry \cite{B}) works efficiently on 
finite dimensional complete Riemannian manifolds.
In contrast to finite dimensions, 
we face a big difficulty to define the $\Gamma_{2}$-operator
when we work
in infinite dimensional frameworks, because it is not so easy
to check the existence of a suitable
core which is not only a ring
but also stable under the operations both of the diffusion semigroup and its generator.
Hence,
we cannot apply this method directly to the infinite dimensional model 
in the present paper.

On the other hand, we have the coupling estimates
(\ref{flow}) and (\ref{flow-2}) which are implied by the strong uniqueness 
of the solution to SPDE (\ref{GL}). By making use of them, 
we can apply the stochastic approach presented 
in \cite{kawa-POTA, Kawa-LSI}.

First, we give the following gradient estimate
for $\{P_{t} \}_{t\geq 0}$.
\begin{pr}[Gradient estimate] \label{thm-GE}
For any $F\in {\cal D(E)}$, we have 
the 
following 
gradient estimate 
\begin{equation} \Vert D(P_{t}F)(w) \Vert_{H}
\leq  
e^{-\frac{K_{1}t}{2}}
P_{t}\big ( \Vert DF \Vert_{H} \big)(w), 
\quad \mu \mbox{-a.e.}~w\in E,~t>0.
\label{GE} 
\end{equation}
\end{pr} 
{\bf Proof:}~The proof is done in the same manner as 
the proof of \cite[Proposition 2.4]{kawa-POTA} together 
with the coupling estimate (\ref{flow-2}). 
So we omit it here.
\qed
\\

Now, we are in a position to state logarithmic Sobolev inequalities.
\begin{tm} 
[Log-Sobolev inequalities]
\label{thm-LSI}
{\rm (1)}~For $F\in {\cal D(E)}$,
we have 
the following heat kernel 
log-Sobolev inequality 
\begin{eqnarray}
& &
\hspace{-25mm}
P_{t}(F^{2}\log F^{2})(w)-P_{t}(F^{2})(w)\log P_{t}(F^{2})(w)
\nonumber \\
& &
\leq \frac{2(1-e^{-K_{1}t})}{K_{1}}P_{t}(\Vert DF \Vert_{H}^{2})(w),
\quad \mu \mbox{-a.e.}~w\in E,~t>0.
\label{Heat-LSI-ineq}
\end{eqnarray}
%
{\rm (2)}~If $K_{1}>0$, that is, $U$ is strictly convex, 
then the following log-Sobolev inequality
\begin{equation}
\int_{E}
F(w)^{2} \log \Big( \frac{F(w)^{2}}{\Vert
F \Vert_{L^{2}(\mu)}^{2}} \Big)\mu(dw)
\leq  \frac{2}{K_{1}}\int_{E} \Vert D_{H}F(w) \Vert_{H}^{2}\mu(dw), \quad F\in 
{\cal D(E)}
\label{LSI-ineq}
\end{equation}
holds.
Consequently, 
we have the spectral gap estimate 
$ \inf \big ( \sigma(-{\cal L}_{\mu})\setminus \{0\} \big ) \geq \frac{K_{1}}{2}$.
\end{tm}
{\bf Proof:}~We first sketch the proof of (1). We refer to \cite{kawa-POTA, Kawa-LSI} for all technical
details. We may assume $F\in {\cal FC}_{b}^{\infty}$, i.e., $F(w)=f(\langle w, \varphi_{1} \rangle,
\ldots, \langle w, \varphi_{n} \rangle)$, where $\{ \varphi_{i} \}_{i=1}^{n} \subset C^{\infty}_{0}
(\mathbb R, \mathbb R^{d})$. Note that $P_{t}F$ can be extended to a function in $C_{b}(E)$ by using
the coupling estimate (\ref{flow}), and the fact that ${\rm{supp}}(\mu)=E$.
We fix $\delta>0$, and introduce a function $G:[0,t] \to L^{1}(\mu)$ by
$$ G(s):=P_{t-s} \big \{ (P_{s}(F^{2})+\delta) \log (P_{s}(F^{2})+\delta) \big \}(\cdot), \quad 0\leq s \leq t.
$$
Then $G$ is differentiable with respect to $s$ and 
\begin{equation}
{\dot G}(s)=-\frac{1}{2} P_{t-s} \Big \{ \frac{ \Vert DP_{s}(F^{2}) \Vert^{2}_{H}}
{P_{s}(F^{2})+\delta } \Big \} (\cdot), \quad 0<s<t.
\label{G-diff}
\end{equation}
On the other hand, Proposition \ref{thm-GE} and Schwarz's inequality imply
\begin{equation}
\Vert DP_{s}(F^{2}) \Vert^{2}_{H} \leq 4e^{-K_{1}s} P_{s}(F^{2}) \cdot P_{s} 
\big( \Vert DF \Vert_{H}^{2} \big).
\label{F2-GE}
\end{equation}
By combining (\ref{G-diff}) with (\ref{F2-GE}), we have
$$ \dot{G}(s) \geq -2e^{-K_{1}s} P_{t-s} \big \{ P_{s} \big( \Vert DF \Vert_{H}^{2} \big) \big \}
=-2e^{-K_{1}s} P_{t} \big( \Vert DF \Vert_{H}^{2} \big).
$$
This imply the heat kernel logarithmic Sobolev inequality (\ref{Heat-LSI-ineq}) 
by first integrating over $s$ from $0$ to $t$ and then
by letting $\delta \searrow 0$. 

Next, we prove (2).
By noting that the Gibbs measure
$\mu$ is the invariant measure for
our stochastic dynamics $\mathbb M$, 
we have the following 
estimate for $w\in E\setminus S$ and $t\geq 0$:
\begin{eqnarray}
\big \vert P_{t}F(w)-{\mathbb E}^{\mu}[F] \big \vert 
&\leq & \int_{E} {\mathbb E} \big [ \vert F(X_{t}^{w})-F(X_{t}^{\widetilde w})
\vert \big ]\mu(d\widetilde w) \nonumber \\
&\leq & \Vert \nabla f \Vert_{\infty} \big (\sum_{i=1}^{n} \Vert \varphi_{i} \Vert_{E^{*}}^{2}
\big )^{1/2}
e^{( \frac{-K_{1}+2 {r}^{2}}{2}t)}
\int \Vert w-{\widetilde w} \Vert_{E}~ \mu(d{\widetilde w}) \nonumber \\
&\leq & {\sqrt 2}
\Vert \nabla f \Vert_{\infty} \big (\sum_{i=1}^{n} \Vert \varphi_{i} \Vert_{E^{*}}^{2}
\big )^{1/2}
e^{( \frac{-K_{1}+2 {r}^{2}}{2}t)}
\Big \{ \Vert w \Vert_{E}^{2}+
\int_{E} \Vert {\widetilde w} \Vert_{E}^{2}\mu(d{\widetilde w}) \Big \}^{1/2},
\nonumber 
\label{ergord}
\\
\end{eqnarray}
where we used (\ref{flow}) for the second line.
Since $r>0$ satisfies $2r^{2}<K_{1}$, (\ref{ergord}) implies the following
ergodic property of $\{P_{t}\}_{t\geq 0}$:
\begin{equation}
\lim_{t\to \infty} P_{t}F(w)={\mathbb E}^{\mu}[F], \quad w\in E\setminus S,
\label{22-ergord}
\end{equation}

Finally, we have the desired logarithmic Sobolev
inequality (\ref{LSI-ineq}) by letting $t\to \infty$ on both sides of (\ref{Heat-LSI-ineq})
and using (\ref{22-ergord}).
This completes the proof of (2).
\qed
\begin{re}
The logarithmic Sobolev inequality {\rm{(\ref{LSI-ineq})}}
holds with $K_{1}\geq m^{2}$ in the case of
exp$(\phi)_{1}$-quantum fields. 
\end{re}
\begin{re}
We mention that 
many other
functional inequalities including the dimension 
free parabolic Harnack inequality (cf. {\rm{\cite{kawa-POTA}}})
and the Littlewood--Paley--Stein inequality 
(cf. Kawabi--Miyokawa {\rm{\cite{KM}}}) 
for our infinite dimensional model 
can be obtained from
the gradient estimate {\rm{(\ref{GE})}}. 
In particular, it is a fundamental and important problem in harmonic analysis
and potential theory
to ask for boundedness of
the Riesz transform $R_{\alpha}({\cal L}_{p}):=D_{H}(\alpha -{\cal L}_{p})^{-1/2}$
on $L^{p}(\mu)$ for all $p>1$ and some $\alpha >0$, where ${\cal L}_{p}$ is the 
extension of $({\cal L}_{0}, {\cal FC}^{\infty}_{b})$ in $L^{p}(\mu)$, because boundedness
of $R_{\alpha}({\cal L}_{p})$ yields the Meyer equivalence of first order Sobolev norms.
In {\rm{\cite{shige-OJM}}}, Shigekawa studied this problem in a general framework
assuming the intertwining
property of the diffusion semigroup $\{P_{t} \}_{t\geq 0}$ and 
another semigroup $\{ \overrightarrow {P}_{t}
\}_{t\geq 0}$ acting on vector-valued functions.
We note that essential self-adjointness of $({\cal L}_{0}, {\cal FC}^{\infty}_{b})$ 
as obtained in Theorem {\rm{\ref{ES}}}
plays a crucial role to prove this property for our model. 
(See e.g., Shigekawa {\rm{\cite{Shige2}}} and
Kawabi {\rm{\cite{RIMS}}}.) We will
discuss boundedness of the Riesz transform by making use of
the Littlewood--Paley--Stein inequality and 
this intertwining property in a forthcoming paper.
\end{re}
\section{
Appendix:~
Another Approach to the Log-Sobolev Inequality (\ref{LSI-ineq}) }
In this section, 
we give another approach to the log-Sobolev inequality (\ref{LSI-ineq}).
First, we prepare the following lemma taken from
Arai--Hirokawa \cite[Lemma 4.9]{Hirokawa}:
\begin{lm}
\label{Hirokawa}
Let 
$\{T_{n} \}_{n=1}^{\infty}$
and $T$ be self-adjoint operators on a Hilbert space
${\cal H}$ having a common core ${\cal D}$ such that, for all $\psi \in {\cal D}$,
$T_{n}\psi \to T\psi$ as $n \to \infty$. Let $\psi_{n}$ be a normalized eigenvectors
of $T_{n}$ with eigenvalue $E_{n}: T_{n}\psi_{n}=E_{n}\psi_{n}$. Assume that
$E:=\lim_{n\to \infty}E_{n}$ exists and that the weak limit ${\mbox{w-}}\lim_{n \to \infty}\psi_{n}=\psi$ 
also exists and one has $\psi \neq 0$. Then $\psi$ is an eigenvector of $T$ with eigenvalue $E$. In particular,
if $\psi_{n}$ is a ground state of $T_{n}$, then $\psi$ is a ground state of $T$.  
\end{lm}
%
\begin{lm} \label{Conv-GS}
Let $U_{N}(z):=\frac{K_{1}}{2}\vert z \vert^{2}+V_{1/N}(z), N=1,2, \ldots$, be potential functions, 
where $V_{1/N}$ is the Moreau--Yosida approximation of $V$.
We consider the Schr\"odinger operator $H_{U_{N}}=-\frac{1}{2}\Delta_{z}+U_{N}$ on 
$L^{2}({\mathbb R}^{d}, {\mathbb R})$, and denote by
$(\lambda_{0})_{N}$ and $\Omega_{N}$ the minimal eigenvalue and the (normalized) ground state of 
$H_{U_{N}}$, respectively. Then the following properties hold under the assumption $K_{1}>0$:
\\
{\rm{(1)}} $(\lambda_{0})_{N} \nearrow \lambda_{0}$ as $N \to \infty$.
\\
{\rm{(2)}} There exists a sub-sequence $\{N(k)\}_{k=1}^{\infty}$ of $N \to \infty$
such that
$\Vert \Omega_{N(k)}-\Omega \Vert_{L^{2}({\mathbb R}^{d}, {\mathbb R})} \to 0$ as $k \to \infty$.
\end{lm}
%
{\bf Proof.}~(1)~Since $U_{N} \nearrow U$ as $N\to \infty$, we have
$(\lambda_{0})_{1}\leq (\lambda_{0})_{2} \leq \cdots \leq \lambda_{0}$.
Moreover, recalling (\ref{falloff})
and taking into account the estimate 
$(U_{N})_{\frac{1}{2}\vert z\vert }(z)\geq \frac{K_{1}}{8}\vert z \vert^{2}, z\in \mathbb R^{d}$ 
for every $N\in \mathbb N$,
we have the following uniform pointwise upper bound for $\{\Omega_{N} \}_{N=1}^{\infty}$:
\begin{equation}
0<\Omega_{N}(z) \leq D_{1} \exp (-\frac{D_{2} K_{1}^{1/2}}{2\sqrt 2} \vert z \vert^{2}), \quad z\in \mathbb R^{d}. 
\label{Omega-N-upper}
\end{equation}

On the other hand, the variational characterization of the minimal eigenvalue and 
the ground state implies
\begin{eqnarray}
(\lambda_{0})_{N}&=& \big({\Omega}_{N},H_{U_{N}}\Omega_{N} \big)_{L^{2}({\mathbb R}^{d}, {\mathbb R})}
\nonumber \\
&=& \big({\Omega}_{N},H_{U}\Omega_{N} \big)_{L^{2}({\mathbb R}^{d}, {\mathbb R})}-
\big({\Omega}_{N},(V-V_{N})\Omega_{N} \big)_{L^{2}({\mathbb R}^{d}, {\mathbb R})}
\nonumber \\
&\geq & \lambda_{0}-
\big({\Omega}_{N},(V-V_{N})\Omega_{N} \big)_{L^{2}({\mathbb R}^{d}, {\mathbb R})},
\label{3-2-1}
\end{eqnarray}
and by Lebesgue's monotone convergence theorem, we also have
\begin{eqnarray}
0&\leq & \lim_{N \to \infty} 
\big({\Omega}_{N},(V-V_{N})\Omega_{N} \big)_{L^{2}({\mathbb R}^{d}, {\mathbb R})}
\nonumber \\
&\leq & 
D_{1}^{2}
\lim_{N \to \infty}  \int_{\mathbb R^{d}} (V(z)-V_{N}(z)) 
\exp (-\frac{D_{2}K_{1}^{1/2}}{\sqrt 2} \vert z \vert^{2})
dz
\nonumber \\
&=&
D_{1}^{2}
\int_{\mathbb R^{d}} \lim_{N \to \infty} (V(z)-V_{N}(z)) 
\exp (-\frac{D_{2}K_{1}^{1/2}}{\sqrt 2} \vert z \vert^{2}) dz=0,
\label{3-2-2}
\end{eqnarray}
where we used (\ref{Omega-N-upper}) for the second line.

Hence by combining (\ref{3-2-1}) with (\ref{3-2-2}), we have 
$\lim_{N \to \infty} (\lambda_{0})_{N} \geq \lambda_{0}$, which completes the proof of (1).
\vspace{1mm} \\
(2)
We take $C^{\infty}_{0}(\mathbb R^{d}, \mathbb R)$ as a common core of the Schr\"odinger
operators $\{H_{U_{N}} \}_{N=1}^{\infty}$ and $H_{U}$ (cf. \cite[Theorem X.28]{rs}), 
and by Lebesgue's monotone convergence theorem, 
we can easily see that
$H_{U_{N}}\psi \to H_{U}\psi$ for all $\psi \in C^{\infty}_{0}(\mathbb R^{d}, \mathbb R)$
as $N \to \infty$. Since $\Vert \Omega_{N} \Vert_{L^{2}({\mathbb R}^{d}, {\mathbb R})}=1$
for all $N\in \mathbb N$, there exist a sub-sequence $\{N(k) \nearrow \infty \}$ and 
a function $\psi \in L^{2}({\mathbb R}^{d}, {\mathbb R})$ such that 
$\Omega_{N(k)} \to \psi$ weakly as $k \to \infty$.
On the other hand, by \cite[Theorem 25.16]{simon}, 
there exist some positive constants $D_{3}, D_{4}$ independent of $N$ such that
\begin{equation}
\Omega_{N}(z) \geq D_{3} \exp \big (-D_{4} \vert z \vert {{U_{N}^{(\infty)}(z)^{1/2}}} \big),
\quad z\in \mathbb R^{d},
\label{Omega-N-lower}
\end{equation}
where $U_{N}^{(\infty)}(z):=
\inf \{ U_{N}(y)\vert~\vert y \vert 
\leq 3 \vert z\vert \}$. Recalling condition {\bf (U3)}, we see that
\begin{eqnarray}
U_{N}^{(\infty)}(z) & \leq &
\inf \{ U(y)\vert~\vert y \vert 
\leq
3 \vert z\vert \}
\nonumber \\
& \leq &
\vert U(0) \vert +3K_{3} \vert z \vert
\exp (3^{\beta} K_{4} \vert z \vert^{\beta} ),
\quad z\in \mathbb R^{d}.
\label{Omega-N-lower-2}
\end{eqnarray}
Then combining (\ref{Omega-N-lower}) with (\ref{Omega-N-lower-2}), we deduce that
\begin{eqnarray}
\Omega_{N}(z) & \geq &
D_{3}
\exp \Big \{ -D_{4}\vert z \vert
{\sqrt{\vert U(0) \vert +3K_{3} \vert z \vert
\exp (3^{\beta} K_{4} \vert z \vert^{\beta} )}}
\Big \} 
\nonumber \\
&=:&\Psi (z), 
\quad z\in \mathbb R^{d}, 
\label{Omega-N-lower-3}
\end{eqnarray}
and hence the uniform pointwise lower estimate (\ref{Omega-N-lower-2}) implies that
$$\lim_{k \to \infty}(\Omega_{N(k)},\Psi)_{L^{2}({\mathbb R}^{d}, {\mathbb R})} 
\geq \Vert \Psi \Vert^{2}_{L^{2}({\mathbb R}^{d}, {\mathbb R})}
>0$$ 
and we now see that $\psi \neq 0$ holds. 

Now by item (1) and Lemma \ref{Hirokawa}, it follows that
$\psi$ is a 
ground state of $H_{U}$. However, since we already know the uniqueness of the ground state of $H_{U}$, 
$\Omega_{N(k)} \to \psi=\Omega$ weakly as $k\to \infty$.
Moreover since $$\lim_{k \to \infty}
\Vert \Omega_{N(k)} \Vert_{L^{2}({\mathbb R}^{d}, {\mathbb R})}
=\Vert \Omega \Vert_{L^{2}({\mathbb R}^{d}, {\mathbb R})}=1,$$ 
we conclude that
$\lim_{k \to \infty} \Vert \Omega_{N(k)}-\Omega \Vert_{L^{2}({\mathbb R}^{d}, {\mathbb R})}=0$.
%
\qed
\vspace{2mm} \\
{\bf Proof of the Log-Sobolev Inequality (\ref{LSI-ineq}):} 
By the same procedure as in Section 2,
we can construct a Gibbs measure $\mu_{N}$ with $\mu_{N}({\cal C})=1$ 
if we replace $U$ by $U_{N}$.
As we have seen in the proof of Theorem \ref{ES}, 
${\widetilde \nabla} U_{1/N}(z)=K_{1}z+\partial_{0} (V_{1/N})(z), z\in {\mathbb R}^{d}$, is Lipschitz continuous. 
Thus, we can apply
\cite[Theorem 1.2]{Kawa-LSI}, and we see
that the following
logarithmic Sobolev inequality holds for each $\mu_{N}$:
\begin{equation}
\int_{E}
F(w)^{2} \log \Big( \frac{F(w)^{2}}{\Vert
F \Vert_{L^{2}(\mu_{N})}^{2}} \Big)\mu_{N}(dw)
\leq  \frac{2}{K_{1}}\int_{E} \Vert D_{H}F(w) \Vert_{H}^{2}\mu_{N}(dw), \quad F\in {\cal FC}^{\infty}_{b}.
\label{LSI-ineq-N}
\end{equation}

Next, we aim to prove the tightness of the family of probability measures $\{\mu_{N} \}_{N=1}^{\infty}$ on ${\cal C}$. 
Due to \cite[Lemma 5.4]{Iwa2},
it suffices to verify the following two conditions: 
\begin{eqnarray}
& & \hspace{-8mm}
{\mbox{ (i)~There exists a constant }} \gamma>0 {\mbox{ such that }}
\sup_{N\in \mathbb N} \int_{\cal C} \vert w(0) \vert^{\gamma} \mu_{N}(dw)<\infty.
\nonumber \\
& & \hspace{-8mm}
{\mbox{ (ii)~For each }}r>0, {\mbox{there exist constants }}p, q, M>0 
{\mbox{ independent of }}N {\mbox{ such that }} 
\nonumber \\
& & 
\int_{\cal C} \vert w(x_{1})-w(x_{2}) \vert^{p} \mu_{N}(dw)\leq M \vert x_{1}-x_{2} \vert^{2+q} \rho_{r}(x_{1}) 
{\mbox{ for }}x_{1},x_{2} \in \mathbb R {\mbox{ with }} \vert x_{1}-x_{2} \vert \leq 1.  
\nonumber
\end{eqnarray}
By combining the translation invariance of $\mu_{N}$ with estimate (\ref{Omega-N-upper}), we have
\begin{equation}
\int_{\cal C} \vert w(0) \vert^{2} \mu_{N}(dw)
=\int_{\mathbb R^{d}} \vert z \vert^{2} \Omega_{N}(z)^{2}dz
\leq D_{1}^{2}\int_{\mathbb R^{d}} \vert z \vert^{2} \exp (-\frac{D_{2}K_{1}^{1/2}}{\sqrt 2} \vert z \vert^{2})dz.
\nonumber
\end{equation}
Hence we have shown that condition (i) holds with $\gamma=2$.
Besides, in a similar way to \cite{Iwa1}, 
we see that
\begin{eqnarray}
& &
\hspace{-10mm}
\int_{\cal C} \vert w(x_{1})-w(x_{2}) \vert^{2m} \mu_{N}(dw)
\nonumber \\
& &
\leq 
\exp \big \{ \big ( (\lambda_{0})_{N}-\inf_{z\in \mathbb R^{d}} U_{N}(z) \big) 
\vert x_{1}-x_{2} \vert \big \} \big( \sup_{z \in \mathbb R^{d}} \Omega_{N}(z) \big)
\nonumber \\
& &
\quad \times
\int_{\mathbb R^{d}} \Omega_{N}(z) dz \cdot (2m-1)!! \cdot \vert x_{1}-x_{2} \vert^{m}
\nonumber \\
& &
\leq
\exp \big \{ \big ( \lambda_{0}-\inf_{z\in \mathbb R^{d}} U_{1}(z) \big) \vert x_{1}-x_{2} \vert \big \}
D_{1}^{2}
%
%
\Big (\frac{ {\sqrt 2}\pi}{D_{2}K_{1}^{1/2}} \Big)^{d/2}
(2m-1)!! \cdot \vert x_{1}-x_{2} \vert^{m}
\nonumber
\end{eqnarray}
for every $m\in \mathbb N$, where $(2m-1)!!:=\prod_{k=1}^{m} (2k-1)$ and
we used Lemma \ref{Conv-GS} and (\ref{Omega-N-upper}) for the third line.
Hence we can find a positive constant $C$ independent of $N$ such that
\begin{equation}
\int_{\cal C} \vert w(x_{1})-w(x_{2}) \vert^{2m} \mu_{N}(dw) \leq C \vert x_{1}-x_{2} \vert^{m},
\quad 
{\mbox{for }}x_{1},x_{2} \in \mathbb R {\mbox{ with }} \vert x_{1}-x_{2} \vert \leq 1, 
\nonumber
\end{equation}
and hence we have proven condition (ii). 

Thus we can find a sub-sequence $\{N(j) \nearrow \infty \}$
such that $\mu_{N(j)}$ converges to some probability measure $\mu_{*}$ weakly on ${\cal C}$. 
On the other hand, by virtue of the Feynman--Kac formula, we have
\begin{equation}
\lim_{N \to \infty} \Vert e^{-tH_{U_{N}}}\psi -e^{-tH_{U}}\psi
\Vert_{L^{2}(\mathbb R^{d}, \mathbb R)}=0, \quad \psi \in L^{2}(\mathbb R^{d}, \mathbb R).
\label{FKsemigroup-conv}
\end{equation}
Then by putting (\ref{FKsemigroup-conv}) and Lemma \ref{Conv-GS} into (\ref{Construction-Gibbs}),
we see that there exists
a sub-sequence $\{N(k) \nearrow \infty \}$ of $\{N(j)\}$ such that
$\lim_{k \to \infty} \mu_{N(k)}(A)=\mu(A)$ for each cylinder set $A\in {\cal B}_{[T_{1},T_{2}]}, T_{1}<T_{2}$.
Hence we obtain $\mu_{*}=\mu$. 

Finally, since $F\in {\cal FC}_{b}^{\infty}$ 
can be regarded as an element of $C_{b}({\cal C})$ in a natural way, 
we can take the limit $k \to \infty$ on both sides of 
(\ref{LSI-ineq-N}). This implies the desired inequality (\ref{LSI-ineq}).
\qed 
\vspace{2mm} \\
{\bf Acknowledgment.}
The authors are grateful to Masao Hirokawa for useful discussions on the paper {\cite{Hirokawa}},
and to Volker Betz and Martin Hairer for providing helpful comments on Remark \ref{Betz-Hairer}.
They were partially supported by the DFG--JSPS joint research 
project \lq \lq Dirichlet Forms, Stochastic Analysis and Interacting Systems" (2007--2008),
and CRC 701, as well as by the project NEST of the Provincia Autonoma di Trento,
at 
University of Trento 
and by HCM at
University of Bonn.
The second named author was also partially supported by Grant-in-Aid for Young Scientists
(Start-up) (18840034) and 
(B) (20740076) from MEXT. This work was completed while the authors were visiting Isaac Newton Institute
for Mathematical Sciences at University of Cambridge. 
They would like to thank the institute for its warm hospitality.


\begin{thebibliography}{99}
\bibitem{A}
S. Albeverio:
{\it{Theory of Dirichlet forms and applications}},
Lectures on probability theory and statistics (Saint-Flour, 2000), 
Lecture Notes in Mathematics, {\bf 1816}, Springer, Berlin, 2003, pp. 1--106, 
%
\bibitem{AFS}
S. Albeverio, F. Flandoli and Y. Sinai:
SPDE in hydrodynamic: recent progress and prospects. 
Lecture Notes in Mathematics, {\bf 1942}. 
Springer-Verlag, Berlin; Fondazione C.I.M.E., Florence, 2008. 
%
\bibitem{AGHk}
S. Albeverio, G. Gallavotti and R. H\o egh-Krohn:
{\it{Some results for the exponential interaction in two or more dimensions}},
Comm. Math. Phys. {\bf 70} (1979), 
pp. 187--192. 
%
\bibitem{AHR}
A. Albeverio, Z. Haba and F. Russo:
{\it{A two-dimensional semilinear heat equation perturbed by (Gaussian)
white noise}},
Probab. Theory Related Fields {\bf{121}} (2001), 
pp. 319--366.
%
\bibitem{AHk73}
S. Albeverio and R. H\o egh-Krohn:
{\it{Uniqueness of the physical vacuum and the Wightman functions in the
infinite volume limit for some non-polynomial interactions}}, 
Commun. Math. Phys. {\bf{30}} (1973), pp. 171--200. 
%
\bibitem{AHk}
S. Albeverio and R. H\o egh-Krohn:
{\it{The Wightman axioms and the mass gap for strong interactions 
of exponential type in two-dimensional space-time}}, 
J. Funct. Anal. {\bf{16}} (1974), pp. 39--82. 
%
%
%
\bibitem{AKKR-book}
S. Albeverio, Y. Kondratiev, Y. Kozitsky and M. R\"ockner:
Statistical Mechanics of Quantum Lattice Systems: A Path Integral Approach,
EMS Tracts in Mathematics {\bf{8}}, European Mathematical Society, 2009.
%
\bibitem{ALZ}
S. Albeverio, S. Liang and B. Zegarli\'nski:
{\it{Remark on the integration by parts formula for the $\phi^4_3$-quantum field model}},
Infin. Dimens. Anal. Quantum Probab. Relat. Top. {\bf 9} (2006), 
pp. 149--154. 
%
\bibitem{ARS}
S. Albeverio, H. R\"ockle and V. Steblovskaya:
{\it{Asymptotic expansions for Ornstein--Uhlenbeck semigroups perturbed by
potentials over Banach spaces}}, 
Stochastics and Stochastics Reports, {\bf 69}, (2000), pp. 195--238. 
%
\bibitem{AR1}
S. Albeverio and M. R\"ockner: 
{\it{New developments in theory and application
of Dirichlet forms}}, in ``Stochastic Processes, Physics and Geometry, 
(Ascona and Locarno, 1988)"
(S. Albeverio, G. Casati, U. Cattaneo, D. Merlini and R. Moresi eds.), 
World Sci. Publ., Teaneck, NJ, 1990, pp. 27--76.
%
\bibitem{AR90}
S. Albeverio and M. R\"ockner: 
{\it{Classical Dirichlet forms on topological vector spaces---
closability and a Cameron-Martin formula}},
J. Funct. Anal. {\bf{88}} (1990), pp. 395--436. 
%
\bibitem{AR}S. Albeverio and M. R\"ockner:
{\it{Stochastic differential equations in infinite dimensions: Solutions
via Dirichlet forms}}, 
Probab. Theory Relat. Fields {\bf 89} (1991) pp. 347--386.
%
%
\bibitem{Cornell}
S. Albeverio and M. R\"ockner:
{\it{Dirichlet form methods for uniqueness of martingale problems 
and applications}}, in
``Stochastic analysis, (Ithaca, NY, 1993)"
(M. Cranston and Mark A. Pinsky eds.),
Proc. Sympos. Pure Math., {\bf{57}}, Amer. Math. Soc., 
Providence, RI, 1995, pp. 513--528.
%
%
\bibitem{Hirokawa}
A. Arai and M. Hirokawa:
{\it{On the existence and uniqueness of ground states of a generalized spin-boson
model}},
J. Funct. Anal. {\bf{151}} (1997), pp. 455--503. 
%
\bibitem{B}D. Bakry: 
{\it{On Sobolev and logarithmic Sobolev inequalities
for Markov semigroups}}, in ``New Trends in Stochastic Analysis''
(K. D. Elworthy, S. Kusuoka and I. Shigekawa eds.),  
World Sci. Publishing, River Edge, NJ (1997), pp. 43--75. 
%
%
\bibitem{D}G. Da Prato: 
{\it{Transition semigroups corresponding to Lipschitz dissipative systems}},
Discrete and Continuous Dynamical Systems {\bf{10}} (2004), pp. 177--192.
%
%
\bibitem{DR}G. Da Prato and M. R\"ockner:
{\it{Singular dissipative stochastic equations in Hilbert spaces}}, 
Probab. Theory Relat. Fields
{\bf{124}} (2002), pp. 261--303. 
%
\bibitem{DT}G. Da Prato and L. Tubaro:
{\it{Some results about dissipativity of Kolmogorov operators}}, Czechoslovak
Mathematical Journal
{\bf{51}} (2001), pp. 685--699. 
%
\bibitem{DT2}
G. Da Prato and L. Tubaro:
{\it{Self-adjointness of some infinite-dimensional elliptic operators and application 
to stochastic quantization}},
Probab. Theory Related Fields 
{\bf{118}} (2000), pp. 131--145.
%
\bibitem{DZ}G. Da Prato and J. Zabczyk: Stochastic Equations
in Infinite Dimensions, Encyclopedia of Mathematics 
and its Applications, Cambridge Univ. Press, Cambridge, UK, 1992. 
%
\bibitem{Eb}
A. Eberle: Uniqueness and non-uniqueness of singular diffusion operators,
Lecture Notes in Mathematics {\bf 1718}, 
Springer-Verlag Berlin Heidelberg, 1999.
%
\bibitem{Fro}
J. Fr\"ohlich:
Nonperturbative quantum field theory, Mathematical aspects and applications,
Selected papers,
Advanced Series in Mathematical Physics, {\bf 15},
World Scientific Publishing Co., River Edge, NJ, 1992. 
%
\bibitem{FOT}M. Fukushima, Y. Oshima and M. Takeda:
Dirichlet Forms and Symmetric Markov Processes,  
Walter de Gruyter, 1994. 
%
\bibitem{Funaki}T. Funaki: {\it{Regularity properties for
stochastic partial differential equations of parabolic type}},  
Osaka J. Math.{\bf{28}} (1991), pp. 495--516.
%
\bibitem{HKPS}
T. Hida, H-H. Kuo, J. Potthoff, and L. Streit:
White noise: An infinite-dimensional calculus, 
Mathematics and its Applications, {\bf 253}. 
Kluwer Academic Publishers Group, Dordrecht, 1993. 
%
\bibitem{Hk}
R. H\o egh-Krohn:
{\it{A general class of quantum fields without cut-offs in two space-time dimensions}},
Comm. Math. Phys. {\bf 21} (1971), pp. 244--255. 
%
%
\bibitem{Iwa1}K. Iwata: {\it{Reversible measures of a 
$P(\phi )_{1}$-time evolution}}, 
in ``Prob. Meth. in Math. Phys.~:~Proceedings of
Taniguchi symposium''~(K. It\^o and N. Ikeda eds.), 
Kinokuniya, 1985,
pp. 195--209. 
%
\bibitem{Iwa2}K. Iwata: {\it{An infinite dimensional stochastic
differential equation with state space $C(\mathbb R)$}},
Probab. Theory Relat. Fields {\bf{74}} (1987), pp. 141--159.
%
%
\bibitem{kawa-POTA}
H. Kawabi: {\it{The parabolic Harnack inequality for the
time dependent Ginzburg-Landau type SPDE and its application}},
Potential Analysis {\bf{22}} (2005), pp. 61--84.
%
\bibitem{Kawa-LSI}
H. Kawabi: {\it{A simple proof of the log-Sobolev inequality on a path 
space with Gibbs measures}},
Infin. Dimens. Anal. Quantum Probab. Relat. Top., 
{\bf{9}} (2006), pp. 321--329.
%
\bibitem{RIMS}
H. Kawabi: {\it{Topics on diffusion semigroups on a path space with Gibbs measures}},
in Proceedings of RIMS Workshop on Stochastic Analysis and Applications, 
(M. Fukushima and I. Shigekawa eds.),
RIMS K\^oky\^uroku Bessatsu, {\bf{B6}}, 2008,
pp. 153--165.
%
\bibitem{KM}
H. Kawabi and T. Miyokawa: {\it{The Littlewood-Paley-Stein inequality 
for diffusion processes on general metric spaces}},
J. Math. Sci. Univ. Tokyo
{\bf{14}} (2007), pp. 1--30.
%
\bibitem{KR}
H. Kawabi and M. R\"ockner:
{\it{Essential self-adjointness
of Dirichlet operators on a path space with Gibbs measures 
via an SPDE approach}},
J. Funct. Anal.
{\bf{242}} (2007), pp. 486--518.
%
%
\bibitem{Ku}
S. Kusuoka:
{\it{H\o egh-Krohn's model of quantum fields and the absolute continuity of measures}},
in R. H\o egh-Krohn's Memorial Volume
``Ideas and methods in quantum and statistical physics'', Vol. 2, 
(S. Albeverio, J.E. Fenstad, H. Holden and T. Lindstr\o m eds.),
Cambridge Univ. Press, Cambridge, 1992, pp. 405--424.  
%
\bibitem{Lisk-Rock}
V. Liskevich and M. R\"ockner:
{\it{Strong uniqueness for certain infinite-dimensional 
Dirichlet operators and applications to stochastic quantization}},
Ann. Scuola Norm. Sup. Pisa Cl. Sci., Serie IV, {\bf 27} (1998), no. 1, 
pp. 69--91.
%
\bibitem{MR} Z.-M. Ma and M. R\"ockner:   
{\rm{Introduction to the Theory of (Non-Symmetric) Dirichlet Forms}},
Springer-Verlag, Berlin-Heidelberg-New York, 1992.
%
\bibitem{Carlo} 
C. Marinelli and M. R\"ockner:
{\it{On uniqueness of mild solutions for dissipative stochastic evolution equations}},
2010, arXiv:1001.5413.
%
%
%
\bibitem{Ondre} 
M. Ondrej\'at:
{\it{Uniqueness for stochastic evolution equations in Banach spaces}},
Dissertaiones Math. (Rozprawy Mat.) {\bf{426}} (2004), 63 pages.
%
\bibitem{priola} 
E. Priola:
{\it{On a class of Markov type semigroups in spaces of
uniformly continuous and bounded functions}},
Studia Mathematica {\bf{136}} (1999), pp. 271--295.
%
\bibitem{rs} M. Reed and B. Simon:
Methods of modern mathematical physics, Vol. II, IV.
New York: Academic Press, 1975, 1978.
%
%
\bibitem{Ro-JFA}
M. R\"ockner:
{\it{Traces of harmonic functions and a new path space for the free quantum field}},
J. Funct. Anal.
{\bf{79}} (1988), pp. 211--249.
%
%
\bibitem{Ro}
M. R\"ockner:
{\it{$L^p$-analysis of finite and infinite-dimensional diffusion operators}},
in ``Stochastic PDE's and Kolmogorov equations in infinite dimensions"
(G. Da Prato ed.), Lecture Notes in Mathematics, {\bf 1715}. 
Springer-Verlag, Berlin, 1999, 
pp. 65--116.
%
%
\bibitem{shige} I. Shigekawa: {\it{$L^{p}$ contraction semigroups for vector valued functions}},
J. Funct. Anal. {\bf{147}} (1997), pp. 69--108. 
%
\bibitem{shige-OJM}
I. Shigekawa:
{\it{Littlewood-Paley inequality for a diffusion satisfying the 
logarithmic Sobolev inequality and for the Brownian motion on 
a Riemannian manifold with boundary}},
Osaka J. Math. {\bf 39} (2002), pp. 897--930.
%
\bibitem{Shige2}
I. Shigekawa:
{\it{Defective intertwining property and generator domain}},
J. Funct. Anal. {\bf{239}} (2006), pp. 357--374.
%
\bibitem{Show}
R.E. Showalter:
Monotone operators in Banach space and nonlinear partial differential equations. 
Mathematical Surveys and Monographs {\bf 49},
American Mathematical Society, Providence, RI, 1997. 
%
\bibitem{Si}
B. Simon:
The $P(\phi )_{2}$ Euclidean (quantum) field theory. 
Princeton Series in Physics. Princeton University Press, Princeton, N.J., 1974.
%
\bibitem{simon} B. Simon:
Functional integration and quantum physics. New York: Academic Press, 1979.
\end{thebibliography}
\end{document}